\documentclass[11pt]{article}

\usepackage[margin=1in]{geometry}
\usepackage{amsmath,amssymb,amsthm}
\numberwithin{equation}{section}
\usepackage{bm}
\usepackage{mathtools}
\usepackage{enumitem}
\usepackage{authblk}

\newtheorem{theorem}{Theorem}[section]
\newtheorem{lemma}[theorem]{Lemma}
\newtheorem{proposition}[theorem]{Proposition}
\newtheorem{corollary}[theorem]{Corollary}
\newtheorem{remark}[theorem]{Remark}
\newtheorem{definition}[theorem]{Definition}
\newtheorem{assumption}[theorem]{Assumption}

\newcommand{\E}{\mathbb{E}}
\newcommand{\bE}{\breve{\mathbb{E}}}
\newcommand{\be}{\breve{\mathcal{E}}}
\newcommand{\wE}{\widehat{\mathbb{E}}}
\newcommand{\we}{\widehat{\mathcal{E}}}
\newcommand{\R}{\mathbb{R}}
\newcommand{\Z}{\mathbb{Z}}
\newcommand{\B}{\widetilde{B}}
\newcommand{\V}{\mathbb{V}}

\newcommand{\usigma}{\overline\sigma}
\newcommand{\lsigma}{\underline\sigma}
\newcommand{\sles}{(\Omega,\mathcal{H},\wE)}

\title{Invariance principle for long-memory linear processes under sublinear expectations\thanks{This work was supported by grants from the National Key R\&D Program of China (No. 2024YFA1013502), the NSF of China (Grant No. U23A2064), and the Summit Advancement Disciplines of Zhejiang Province (Zhejiang Gongshang University--Statistics).}}

\author[1]{Wang-Yun Gu}
\author[2]{Li-Xin Zhang}

\affil[1]{School of Mathematical Sciences, Zhejiang University, Hangzhou, China\\
  E-mail:\texttt{guwangyun0408@gmail.com}}

\affil[2]{School of Statistics and Data Science, Zhejiang Gongshang University, Hangzhou, China\\
  E-mail:\texttt{stazlx@mail.zjgsu.edu.cn}}

\date{}
\date{}

\begin{document}

\maketitle

\begin{abstract}

We establish an invariance principle for long-memory linear processes with power-law coefficients under sublinear expectations. Under truncated-moment and capacity-tail conditions, we construct these processes on a canonical space and prove weak convergence in continuous path space of their linearly interpolated, normalized partial sums. The limiting sublinear law is represented as an upper expectation over stochastic integrals with a fractional Brownian moving-average kernel and predictable volatility controls, driven by a two-sided classical Brownian motion. When variance uncertainty vanishes, the limit reduces to the classical fractional Brownian limit. The proof combines a multidimensional Rosenthal inequality, coefficient approximation, and the invariance principle for independent and identically distributed random variables to establish finite-dimensional convergence. Tightness follows from second-moment increment bounds. Uniform approximation on compact sets then extends the convergence from bounded Lipschitz to bounded continuous test functions.
\end{abstract}



\section{Introduction}

Limit theorems for linear processes
\begin{equation*}
Y_i=\sum_{j\in\Z}\alpha_j\varepsilon_{i-j}
\end{equation*}
with i.i.d.\ innovations $\varepsilon_i$ are a classical topic in probability theory~\cite{IbragimovLinnik1971}.  For centered innovations with finite, nonzero variance, in the usual non-degenerate short-range case
$\sum_j|\alpha_j|<\infty$ with $\sum_j\alpha_j\ne0$, the partial sums
$\sum_{i=1}^{\lfloor nt\rfloor}Y_i$ satisfy a functional central limit theorem under
$\sqrt n$-normalization \cite{Hannan1979,PeligradUtev2006}.  If instead
\begin{equation*}
j^\beta\alpha_j\to c_+,\qquad j^\beta\alpha_{-j}\to c_-,\qquad
\beta\in(1/2,1),\ c_+^2+c_-^2>0,
\end{equation*}
then $\sum_j|\alpha_j|=\infty$ but $\sum_j\alpha_j^2<\infty$; the partial sums carry
long-range dependence, the variance grows like $n^{2H}$ with
\begin{equation*}
H=\frac32-\beta\in\left(\frac12,1\right),
\end{equation*}
and the classical fractional-Brownian invariance principle yields
\begin{equation*}
n^{-H}\sum_{i=1}^{\lfloor nt\rfloor}Y_i \Longrightarrow \kappa\,\sigma B^H_t,
\end{equation*}
where $B^H$ is standard fractional Brownian motion (fBm) with Hurst parameter $H$, $\sigma^2=E[\varepsilon_0^2]$, and $\kappa$ is a
constant determined by the filter~\cite{Davydov1970, DMP2011}.  Classical long-memory limit theory exploits a range of
tools tied to linear probability, including characteristic-function arguments, spectral and
multiple Wiener--It\^o methods, and Appell-polynomial expansions; see, for example, Taqqu
\cite{Taqqu1975}, Dobrushin and Major \cite{DobrushinMajor1979}, Major \cite{Major1981}, and
Avram and Taqqu \cite{AvramTaqqu1987}.  For linear processes, the later
invariance-principle argument of Dedecker, Merlev\`ede and Peligrad \cite{DMP2011} is based
instead on martingale approximation.

The situation is different under a \emph{sublinear} expectation
$\wE$ in the sense of Peng \cite{Peng2019}.  
Classical characteristic-function arguments do not carry
over directly to the sublinear-expectation setting
\cite{BayraktarMunk2016,HuJiangLiangPeng2023}.
For linear processes with absolutely summable filters
under sublinear expectations, Zhang
\cite{Zhang2021} proves a central limit
theorem for one-sided moving averages with i.i.d.\
innovations. Liu and Zhang \cite{LiuZhang2021} establish
a central limit theorem and an invariance principle
for linear processes generated by i.i.d. innovations, while Guo and Zhang
\cite{GuoZhang2023} obtain a central limit theorem
for $m$-dependent innovations. 
We consider the long-memory regime with square summable filters, establishing
an invariance principle under the normalization
$n^H$, where $H\in(1/2,1)$.

Under the same moment conditions on the innovations as in Zhang’s i.i.d. functional central limit theorem \cite{Zhang2023},
we prove weak convergence of the linearly interpolated,
normalized partial-sum process on $C[0,1]$, with convergence
of sublinear expectations for every bounded continuous
test function.
The limiting sublinear law is represented as an upper
expectation over stochastic integrals with a fractional
Brownian moving-average kernel and predictable volatility
controls, driven by a two-sided classical Brownian motion.

We first establish a multidimensional Rosenthal inequality
and use its second-moment case to construct the linear
process in a complete $L^2$ space on a canonical realization,
without assuming a finite second moment on the original space.
Coefficient approximation and the i.i.d.\ invariance principle
give convergence for weighted sums over bounded intervals;
tail estimates then yield finite-dimensional convergence.
Finally, second-moment increment bounds yield tightness,
while dyadic interpolation and approximation on compact sets
establish weak convergence against all bounded continuous
test functions.

The paper is organized as follows.
Section~\ref{sec:setup} introduces the framework and assumptions.
Section~\ref{sec:main} describes the limiting law and states the
main theorem and its classical specialization.
Section~\ref{sec:rosenthal} develops the multi-dimensional Rosenthal's inequalities used in the well-posedness and 
finite-dimensional argument.
Section~\ref{sec:weak-wellposedness} establishes the well-posedness of the linear process under the truncated-moment assumptions.
Section~\ref{sec:fdc} proves finite-dimensional convergence.
Section~\ref{sec:tightness} establishes tightness and completes the proof
of the main theorem.

\section{Setup and assumptions}\label{sec:setup}
We use the framework of Peng~\cite{Peng2019}. Let $(\Omega,\mathcal{F})$ be a given measurable space and let $\mathcal{H}$ be a linear space of real functions defined on $(\Omega,\mathcal{F})$ such that if $X_1,\cdots,X_n\in\mathcal{H}$ then $\varphi(X_1,\cdots,X_n)\in \mathcal{H}$ for each $\varphi\in C_{l,Lip}(\mathbb{R}^n)$, where $C_{l,Lip}(\mathbb{R}^n)$ denotes the linear space of locally Lipschitz functions $\varphi$ satisfying 
$$|\varphi(\bm{x})-\varphi(\bm{y})|\leq C(1+|\bm{x}|^m+|\bm{y}|^m)|\bm{x}-\bm{y}|,\quad \forall \bm{x},\bm{y}\in\R^n,$$
for some $C>0,m\in\mathbb{N}$ depending on $\varphi$. $\mathcal{H}$ is considered as a space of "random variables". We also denote the space of bounded Lipschitz functions by $C_{b,Lip}(\R^n)$. A sublinear expectation $\wE$ on $\mathcal{H}$ is a function
$\wE:\mathcal H\to[-\infty,\infty]$ satisfying the following four properties: for all $X,Y\in\mathcal{H}$,
\begin{itemize}
    \item[(1)]Monotonicity: If $X\ge Y$, then $\wE[X]\ge \wE[Y]$;
    \item[(2)]Constant preserving: $\wE[c]=c$;
    \item[(3)]Sub-additivity: $\wE[X+Y]\le\wE[X]+\wE[Y]$ whenever  $\wE[X]+\wE[Y]$ is well-defined.
    \item[(4)]Positive homogeneity: $\wE[\lambda X]=\lambda\wE[X],\enspace\forall\lambda\ge0$.
\end{itemize}
The triple $\sles$ is called a sublinear expectation space. Given a sublinear expectation $\wE$, we define its conjugate expectation $\we$ by $\widehat{\mathcal E}[X]:=-\wE[-X],\enspace \forall X\in\mathcal H.$ As discussed in Remark 1.4.5 of \cite{Peng2019}, when we talk about the completion space $\mathcal{H}$ under the norm $\wE[|X|^p]^{1/p}$,  $X_1,\cdots,X_n\in\mathcal{H}$ does not imply in general that $\varphi(X_1,\cdots,X_n)\in \mathcal{H}$ for each $\varphi\in C_{l,Lip}(\mathbb{R}^n)$. Thus we define the notations of distribution and independence on $C_{b,Lip}(\R^n)$.

\begin{definition}
\label{def:independence}
(see Peng~\cite{Peng2019})\\
\begin{itemize}
    \item [(1)] (Identical distribution) Let $\bm{X}_1$ and $\bm{X}_2$ be two $n$-dimensional random vectors defined respectively in sublinear expectation spaces $(\Omega_1, \mathcal{H}_1,\wE_1)$ and $(\Omega_2, \mathcal{H}_2,\wE_2)$. They are called identically distributed, denoted by $\bm{X}_1\overset{d}{=}\bm{X}_2$ if
    \begin{equation*}
        \wE_1[\varphi(\bm{X}_1)]=\wE_2[\varphi(\bm{X}_2)],\quad \forall \varphi\in C_{b,Lip}(\R^n).
    \end{equation*}
    A sequence of random variables $\{X_n;n\in\Z\}$  is said to be identically distributed if $X_i\overset{d}{=}X_1$ for each $i\in\Z$.
    \item[(2)] (Independence) In a sublinear expectation space $\sles$, a random vector $\bm{Y}=(Y_1,\cdots,Y_n)$, $Y_i\in\mathcal{H}$ is said to be independent of another random vector $\bm{X}=(X_1,\cdots,X_m), X_i\in\mathcal{H}$ under $\wE$ if for each test function $\varphi\in C_{b,Lip}(\R^m\times \R^n)$, we have  
    \begin{equation*}
        \wE[\varphi(\bm{X},\bm{Y})]=\wE[\wE[\varphi(\bm{x},\bm{Y})]|_{\bm{x}=\bm{X}}].
    \end{equation*}
    A sequence of random variables $\{X_n;n\in\Z\}$ is said to be independent, if $X_{j_r}$ is independent of $(X_{j_1},\cdots, X_{j_{r-1}})$ whenever $j_1<\cdots < j_r$.
\end{itemize}
\end{definition}

We also recall condition \textup{(CC)}; see
\cite[Lemma~2.5(a)]{Zhang2023SLLN}.
Write $\mathcal H_b:=\{X\in\mathcal H: X\text{ is bounded}\}.$
We say that $\wE$ satisfies \textup{(CC)} if there exists
a nonempty family $\mathcal P$ of probability measures
on $(\Omega,\mathcal F)$ such that
\begin{equation}\label{eq:CC-representation}
 \wE[X]
 =
 \sup_{P\in\mathcal P}E_P[X],
 \qquad X\in\mathcal H_b,
\end{equation}
and $\mathcal P$ is countable-dimensionally weakly
compact.
Precisely, for every sequence $\{X_j;j\ge 1\}$ in
$\mathcal H_b$ and every sequence $\{P_n\}$ in
$\mathcal P$, there exist a subsequence $\{P_{n_k}\}$
and $P\in\mathcal P$ such that
\[
 \lim_{k\to\infty}
 E_{P_{n_k}}[\varphi(X_1,\ldots,X_d)]
 =
 E_P[\varphi(X_1,\ldots,X_d)]
\]
for every $d\ge1$ and
$\varphi\in C_{b,\mathrm{Lip}}(\R^d)$.

Now we introduce the capacity. Let $\mathcal{G}\subset\mathcal{F}$. A function $V:\mathcal{G}\mapsto[0,1]$ is called a capacity if 
\begin{equation*}
    V(\varnothing)=0, \enspace V(\Omega)=1\text{ and } V(A)\le V(B) \; \forall A\subset B,\;A,B\in\mathcal{G}.
\end{equation*}
It is called sub-additive if $V(A\cup B)\le V(A)+V(B)$ for all $A,B\in\mathcal{G}$ with $A\cup B\in\mathcal{G}$.

Let $\sles$ be a sublinear expectation space. Let $(\V,\mathcal{V})$ be a pair of capacities with the properties that
\begin{align}
    \wE[f]\le \V(A)\le\wE[g],\;\;\text{if }f\le I_A\le g, \;f,g\in\mathcal{H}\text{ and }A\in\mathcal{F},\label{eq:propertyV}
\end{align}
$\V \text{ is sub-additive}$ and $\mathcal{V}(A)=1-\V(A^c),\;\forall A\in\mathcal{F}.$
We call $\V$ and $\mathcal{V}$ the upper and the lower capacity, respectively. We say that an event holds quasi-surely if $A\in\mathcal{F}$ with $\mathbb V(A^c)=0$.

When there exists a family of probability measures on $(\Omega,\mathcal{F})$ such that \eqref{eq:CC-representation} holds,
$\V$ can be defined as 
\begin{equation}
    \V(A)=\sup_{P\in\mathcal{P}}P(A).
\end{equation}

If $\V_1$ and $\V_2$ are two capacities with the property \eqref{eq:propertyV}, then for any random variables $X\in\mathcal{H}$, 
\begin{equation}
    \V_1(X\ge x+\varepsilon)\le \V_2(X\ge x)\le \V_1(X\ge x-\varepsilon)\;\; \text{for all }\varepsilon>0 \text{ and }x.
\end{equation}

Moreover, we recall the notations of $G$-normal distribution and $G$-Brownian motion introduced in Peng~\cite{Peng2019}.
\begin{definition}[$G$-normal random variable]
    For $0\le\lsigma^2\le\usigma^2<\infty$, a random variable $\xi$ in a sublinear expectation space $(\widetilde{\Omega},\widetilde{\mathcal{H}},\widetilde{{\E}})$ is called a normally $N(0,[\lsigma^2,\usigma^2])$ distributed random variable (written as $\xi\sim N(0,[\lsigma^2,\usigma^2])$ under $\widetilde{\E}$), if for any $\varphi\in C_{b,Lip}(\R)$, the function $u(t,x)=\widetilde{\E}[\varphi(x+\sqrt{t}\xi)] $ is the unique viscosity solution of the following heat equation:
    $$\partial_tu-G(\partial^2_{xx}u)=0,\;\;\;u(0,x)=\varphi(x), \;\;x\in\R,t\ge 0,$$
    where $G(\alpha)=\frac12(\usigma^2\alpha^+-\lsigma^2\alpha^-)$.
\end{definition}
That $\xi$ is a normally distributed random variable is equivalent to that, if $\xi'$ is an independent copy of $\xi$,
\begin{equation*}
        \widetilde{\E}[\varphi(a\xi+b\xi')]=\widetilde{\E}[\varphi(\sqrt{a^2+b^2}\xi)],\;\;\;\forall \varphi\in C_{b,Lip}(\R), \text{ and }a,b\ge0.
\end{equation*} 
\begin{definition}[$G$-Brownian motion]
    A random process $(B_t)_{t\ge 0}$ in the sublinear expectation space $(\widetilde{\Omega},\widetilde{\mathcal{H}},\widetilde{{\E}})$ is called a $G$-Brownian motion if 
    \begin{itemize}
        \item [(1)] $B_0=0$;
        \item[(2)] for each $0\le t_1\le \cdots\le t_d\le t\le s$ and $\varphi\in C_{b,Lip}(\R^{d+1})$,
        \begin{equation}
            \widetilde{\E}[\varphi(B_{t_1},\cdots,B_{t_d}, B_s-B_t)]=\widetilde{\E}[\widetilde{\E}[\varphi(x_1,\cdots,x_d,\sqrt{s-t}\xi)]|_{(x_1,\cdots,x_d)=(B_{t_1},\cdots, B_{t_d})}],
        \end{equation}
        where $\xi\sim N(0,[\lsigma^2,\usigma^2])$.
    \end{itemize}
Without loss of generality,  we assume that $\widetilde{\Omega}$ is a Polish space and there is a family of probability measures $\mathcal Q$ on $\widetilde{\Omega}$ such that
$$\widetilde{\E}[X]=\sup_{Q\in \mathcal{Q}}E_Q[X], X\in \widetilde{\mathcal{H}}. $$
Furthermore, for any measurable function $f$ on $\widetilde{\Omega}$, define 
$$\widetilde{\E}[f]=\sup_{Q\in \mathcal{Q}}E_Q[f]. $$
\end{definition}

For real numbers $x$ and $y$, denote $x\vee y=\max(x,y),\;x\wedge y=\min(x,y)$, $x^+=\max(x,0)$ and $x^-=\max(0,-x)$. For a random variable $X$, we will truncate it in the form $(-c)\vee X\wedge c$ denoted by $X^{(c)}$, and define $\bE[X]=\lim_{c\to\infty}\wE[X^{(c)}]$ if the limit exists, and $\be[X]=-\bE[-X]$.

Finally, we give the assumptions of the innovations $\{\varepsilon_i\}_{i\in\Z}$ and coefficients $\{\alpha_i\}_{i\in\Z}$ of the long range dependent linear process discussed in the sequel.
\begin{assumption}\label{ass:innovations}
$\{\varepsilon_i\}_{i\in\Z}$ is a sequence of i.i.d.\ random variables under $\wE$.  We assume that
\begin{align}
    &\bE[\varepsilon_1]=\be[\varepsilon_1]=0,\\
&\usigma^2:=\bE[\varepsilon_1^2]<\infty,
\label{eq:moments}
\end{align}
and 
\begin{equation}
    \lim_{x\to\infty}x^2\V(|\varepsilon_1|\geq x)=0,
\label{ass:tight-condition}
\end{equation}
where $\V$ is a capacity with property \eqref{eq:propertyV}.
\end{assumption}

\begin{assumption}\label{ass:coeff}
There are $\beta\in(1/2,1)$ and real constants $c_+,c_-$ with $c_+^2+c_-^2>0$ such that
\begin{equation*}
j^{\beta}\alpha_j\to c_+,\qquad
j^{\beta}\alpha_{-j}\to c_-,\qquad\text{as}\enspace j\to+\infty.
\end{equation*}
In particular $\sum_j|\alpha_j|=\infty$ and $\sum_j\alpha_j^2<\infty$.  We set
\begin{equation*}
H:=\frac32-\beta\in\left(\frac12,1\right).
\end{equation*}
\end{assumption}

Define the linear process and its discrete partial sums by
\begin{equation*}
\begin{aligned}
&Y_i=\sum_{j\in\Z}\alpha_j\varepsilon_{i-j}, \qquad S_m:=\sum_{i=1}^{m}Y_i\quad(m\ge1),\qquad S_0:=0,\\
&X_n(t):=n^{-H}S_{\lfloor nt\rfloor},\quad t\in[0,1].
\end{aligned}
\end{equation*}

Since $Y_i$ may not be an element of $\mathcal {H}$, we understand the series $Y_i$  as an $L_{\mathcal P}^2$ limit on the canonical space constructed
in Section~\ref{sec:weak-wellposedness},  and for any $\varphi\in C_{b,Lip}(\mathbb R^d)$, the sublinear expectation of $\varphi(Y_1,\ldots, Y_d)$ is defined as 
\begin{equation}\label{eq:extestionE:1}
\wE\left[\varphi(Y_1,\ldots,Y_d)\right]= \lim_{N\to \infty}\wE\left[\varphi(Y_1^{(N)},\ldots,Y_d^{(N)})\right], 
\end{equation} 
where 
$Y_i^{(N)}=\sum_{|j|\le N}\alpha_j\epsilon_{i-j}$,
see Corollary~\ref{cor:wp-linear-process} and \eqref{eq:wp-transfer-fdd}. 
Furthermore, for any $\varphi\in C_b(\mathbb R^d)$,
there exists a sequence
$\{\varphi_k\}\subset C_{b,\mathrm{Lip}}(\mathbb R^d)$
such that
\[
\sup_{x\in[-R,R]^d}|\varphi_k(x)-\varphi(x)|
\longrightarrow 0,
\qquad \forall R>0,
\]
and $\sup_{k\ge1}\|\varphi_k\|_\infty \le \|\varphi\|_\infty.$
It can be shown that the limit
$\lim_{k\to\infty}
\widehat{\mathbb E}[\varphi_k(Y_1,\ldots,Y_d)]$
exists and does not depend on the choice of the sequence
$\{\varphi_k\}$ satisfying the above two conditions.
Thus, we can extend  $\wE$ to that on bounded continuous functions  of $(Y_1,\ldots, Y_d)$, $d\ge 1$, by defining 
\begin{equation}\label{eq:extestionE:2}\wE\left[\varphi(Y_1,\ldots, Y_d)\right]=\lim_{k\to \infty}\wE\left[\varphi_k(Y_1,\ldots,Y_d)\right].
\end{equation}

If $ \wE $ has the representation \eqref{eq:CC-representation},
then
$$\wE\left[\varphi(Y_1,\ldots,Y_d)\right]=\sup_{P\in\mathcal P}E_P\left[\varphi(Y_1,\ldots,Y_d)\right], \; \varphi\in C_b(\mathbb R^d). $$

\section{Main result}\label{sec:main}

The linear process and its expectations are understood through the
realization constructed in Section~\ref{sec:weak-wellposedness},
with the original notation retained.
We use the linear interpolation
\begin{equation}\label{eq:defX}
 \widetilde X_n(t)
 :=n^{-H}\bigl[S_{\lfloor nt\rfloor}
 +(nt-\lfloor nt\rfloor)Y_{\lfloor nt\rfloor+1}\bigr],
 \qquad t\in[0,1],
\end{equation}
which coincides with $X_n$ at the grid points.

We first describe the limiting law.
Let $\overline\sigma^2$ be as in
Assumption~\ref{ass:innovations}, and set $\lsigma^2:=\be[\varepsilon^2]$.
The limit exists by monotonicity, and
$0\le\underline\sigma^2\le\overline\sigma^2<\infty$.
Let $W=(W_x)_{x\in\R}$ be a two-sided standard Brownian motion
on a classical probability space
$(\Omega^W,\mathcal F^W,P^W)$, and write
$E:=E_{ P^W}$.
We use the increment filtration
\[
 \mathcal F_x^W
 :=\sigma\{W_v-W_u:u<v\le x\},
 \qquad x\in\R,
\]
with the usual augmentation, and define
\begin{equation}\label{eq:whole-line-controls}
 \mathcal A=\mathcal A_{\R}
 :=\bigl\{\theta:
 \theta\text{ is }(\mathcal F_x^W)\text{-predictable and takes values in }
 [\underline\sigma,\overline\sigma]\bigr\}.
\end{equation}
For $a<b$, let
\[
 \mathcal F_x^{W,a}
 :=\sigma\{W_r-W_a:a\le r\le x\},
 \qquad a\le x\le b,
\]
again with the usual augmentation, and set
\begin{equation}\label{eq:finite-interval-controls}
 \mathcal A_{[a,b]}
 :=\bigl\{\theta:
 \theta\text{ is }(\mathcal F_x^{W,a})\text{-predictable on }[a,b]
 \text{ and takes values in }[\underline\sigma,\overline\sigma]\bigr\}.
\end{equation}
Thus a finite-interval control uses no Brownian noise prior to $a$.

For $t\ge0$, define
\begin{equation}\label{eq:main-limit-kernel}
 f_t(x)
 :=\int_0^t
 \bigl[c_+(u-x)_+^{-\beta}+c_-(x-u)_+^{-\beta}\bigr]\,du,
 \qquad x\in\R,
\end{equation}
where $r_+^{-\beta}$ means $r^{-\beta}$ for $r>0$ and $0$ otherwise.
These are the moving-average kernels of fractional Brownian
motion; see \cite[Eq.~(7.11)]{Samorodnitsky2006LRD}.
In particular, $f_t\in L^2(\mathbb R)$. Set
\[
 \kappa^2:=\int_{\mathbb R}f_1(x)^2\,dx\in(0,\infty).
\]
For $0\le s<t$, translation and scaling give $
 f_t(x)-f_s(x)=f_{t-s}(x-s)
 =(t-s)^{H-\frac12}
 f_1\!\left(\frac{x-s}{t-s}\right).
$
Consequently,
\[
 \|f_t-f_s\|_{L^2(\mathbb R)}^2
 =\kappa^2|t-s|^{2H},
 \qquad s,t\ge0.
\]
For $\theta\in\mathcal A$, set
\[
 B_t^\theta:=\int_{\R}f_t(x)\theta_x\,dW_x,
 \qquad t\ge0,
\]
where the integral is the $L^2( P^W)$-limit of the integrals
over $[-R,R]$ as $R\to\infty$.
By It\^o's isometry,
\[
 E[|B_t^\theta-B_s^\theta|^2]
 \le\overline\sigma^2\|f_t-f_s\|_{L^2(\R)}^2
 =\overline\sigma^2\kappa^2|t-s|^{2H}.
\]
Since $2H>1$, Kolmogorov's continuity theorem gives a continuous
modification, which we use throughout.
For $\Phi\in C_{b}(C[0,1])$, define
\begin{equation}\label{eq:main-limit-path-law}
 \widetilde{\E}[\Phi]
 :=\sup_{\theta\in\mathcal A}
 E\bigl[\Phi((B_t^\theta)_{t\in[0,1]})\bigr].
\end{equation}
We write $\B^H$ for the coordinate process under this sublinear law.

\begin{theorem}\label{thm:main}
Under Assumptions~\ref{ass:innovations} and~\ref{ass:coeff},
$\widetilde X_n$ converges weakly to $\B^H$ in $C[0,1]$,
equipped with the supremum norm: for every
$\Phi\in C_{b}(C[0,1])$,
\begin{equation}\label{eq:main}
 \lim_{n\to\infty}\wE[\Phi(\widetilde X_n)]
 =\widetilde{\E}[\Phi].
\end{equation}
In particular, for every $d\ge1$, $0\le t_1<\cdots<t_d\le1$,
and $\varphi\in C_{b}(\R^d)$,
\begin{equation}\label{eq:limitprocess}
\begin{aligned}
 \lim_{n\to\infty}
 \wE\bigl[\varphi(\widetilde X_n(t_1),\ldots,
                  \widetilde X_n(t_d))\bigr]=
 L[\varphi]
 :=\sup_{\theta\in\mathcal A}
 E\left[\varphi\left(
 \int_{\R}\bm{f}(x)\theta_x\,dW_x
 \right)\right],
\end{aligned}
\end{equation}
where $\bm{f}(x):=(f_{t_1}(x),\ldots,f_{t_d}(x))^\top$.
\end{theorem}
\begin{remark} Notice that for every
$\Phi\in C_{b}(C[0,1])$, 
$\Phi(\widetilde X_n)$ is a bounded continuous function of $S_1,\ldots,S_n$, and so a  bounded continuous function of $Y_1,\ldots,Y_n$. $\wE[\Phi(\widetilde X_n)]$ is well-defined by \eqref{eq:extestionE:1} and \eqref{eq:extestionE:2}.
    
\end{remark}
\begin{corollary}\label{cor:classical}
If $\underline\sigma=\overline\sigma=\sigma$, then
\[
 \widetilde{\E}[\Phi]
 =E[\Phi(\kappa\sigma B^H)],
 \qquad \Phi\in C_{b}(C[0,1]),
\]
where $B^H$ is a standard fractional Brownian motion with Hurst
parameter $H$.
In particular, when $\wE$ is a classical probability expectation,
Theorem~\ref{thm:main} recovers Davydov's invariance principle
\cite{Davydov1970}.
\end{corollary}

\begin{proof}
The control is then constant, $\theta_x=\sigma$.
The resulting process $B^\sigma$ is zero-mean Gaussian, with covariance
\[
 \E[B_t^\sigma B_s^\sigma]
 =\sigma^2\int_{\R}f_t(x)f_s(x)\,dx
 =\frac{\sigma^2\kappa^2}{2}
 \bigl(t^{2H}+s^{2H}-|t-s|^{2H}\bigr).
\]
Thus its law is that of $\kappa\sigma B^H$.
\end{proof}

\section{Multi-dimensional Rosenthal's inequalities}\label{sec:rosenthal}
Zhang~\cite{Zhang2016} established Rosenthal's inequalities
for real-valued random variables under sublinear expectations.
Here we extend these inequalities to random vectors with respect to
the Euclidean norm. We begin with two auxiliary lemmas.

\begin{lemma}[Young's inequality with $\delta$]
\label{lem:Young}
Let $r,s>1$ satisfy $\frac{1}{r}+\frac{1}{s}=1$. Then for all $a,b\ge 0$ and every $\delta>0$,
\begin{equation}
ab\;\le\;\delta\,a^{r}+\frac{\delta^{-s/r}}{s\,r^{\,s/r}}\,b^{s}.
\label{eq:weighted}
\end{equation}
\end{lemma}

\begin{proof}
Fix $\lambda>0$ and replace $(a,b)$ by $(\lambda a,\lambda^{-1}b)$ in Young's inequality:
\begin{equation}
ab=(\lambda a)\bigl(\lambda^{-1}b\bigr)
\;\le\;\frac{\lambda^{r}}{r}\,a^{r}+\frac{\lambda^{-s}}{s}\,b^{s}.
\label{eq:scaled}
\end{equation}
Set $\delta:=\frac{\lambda^{r}}{r}$, i.e.\ $\lambda=(r\delta)^{1/r}$. Then
\begin{equation*}
\lambda^{-s}=\Bigl((r\delta)^{1/r}\Bigr)^{-s}=(r\delta)^{-s/r}=r^{-s/r}\,\delta^{-s/r},
\end{equation*}
and substituting into \eqref{eq:scaled} yields \eqref{eq:weighted}.
\end{proof}

\begin{lemma} [Proposition 2.1 in \cite{Acciaio2013}]
\label{lem:doob}
    Let $\{r_0, r_1, \cdots,r_n\}$ be a non-negative sequence. Define $r_k^*=\max_{0\le j \le k}r_j$, $q=\frac{p}{p-1}$ for some $p>1$. Then 
    \begin{equation}
        (r_n^*)^p\leq q^pr_n^p-pq\sum_{i=1}^n(r_{i-1}^*)^{p-1}(r_i-r_{i-1}).
        \label{eq:max-estimate}
    \end{equation}
\end{lemma}

\begin{theorem}
\label{thm:rosenthal}
    Suppose $\{\bm{X}_n;n\ge 1\}$ is a sequence of independent $d$-dimensional random vectors on $\sles$ and $\bE[\Vert X_i\Vert_2^p]<\infty$ for some $p\ge 2$. For deterministic vectors $\bm{c}_i\in\mathbb{R}^d$, denote $\bm{Y}_i=\bm{X}_i-\bm{c}_i, \bm{T}_0=\bm{0},\bm{T}_k=\sum_{i=1}^k\bm{Y}_i$ and 
    \begin{equation}
        A:=\sum_{i=1}^n\bE[\Vert \bm Y_i\Vert_2^p], \quad B:=\sum_{i=1}^n\bE[\Vert \bm Y_i\Vert_2^2], \quad D:=\sum_{i=1}^n\rho_i(\bm{c}_i)
    \end{equation}
    where $\rho_i(\bm{c}_i)=\max_{\bm\theta\in\Theta_i}\Vert \bm\theta-\bm{c}_i\Vert_2$, $\Theta_i$ is the mean uncertainty set of $\bm{X}_i$ on $\mathbb{R}^d$, i.e.
    \begin{equation}
        \Theta_i=\{\bm{\theta}\in\mathbb{R}^d:\langle\bm{u},\bm{\theta}\rangle\leq\bE[\langle\bm{u},\bm{X}_i\rangle], \enspace\forall \bm{u}\in \mathbb{R}^d, \|\bm u\|_2\le 1\}.
    \end{equation}
    Then there exists a constant $C_p>0$ such that 
    \begin{equation}
    \label{eq:Rosenthal}
        \bE\bigg[\max_{0\le k \le n}\Vert \bm{T}_k\Vert_2^p\bigg]\leq C_p(A+B^{p/2}+D^p).
    \end{equation}
\end{theorem}
\begin{proof} First, we assume that the random vectors $\bm X_n$ are bounded. Notice that for a bounded $\bm X$, we have $\bE[\bm X]=\wE[\bm X]$. 
    Denote $\rho_i=\rho_i(\bm{c}_i)$ and note
    \begin{equation*}
        \wE[\langle \bm{u}, \bm{Y}_i\rangle]=\max_{\bm{\theta}\in\Theta_i}\langle\bm{u},\bm\theta-\bm{c}_i\rangle\leq \Vert\bm{u}\Vert_2\max_{\bm\theta\in\Theta_i}\Vert \bm\theta-\bm{c}_i\Vert_2=\rho_i\Vert \bm{u}\Vert_2.
    \end{equation*}

Since the random variables are bounded, truncating the arguments of the locally Lipschitz test functions used below allows us to apply Definition~\ref{def:independence}.    
For $d$-dimensional function $\bm{H}=\bm{H}(\bm{Y_1,\cdots,\bm{Y}_{i-1}})\in\mathcal{H}$, it follows from independence that
    \begin{align}
    \label{eq:mean-estimate}
        \wE[\pm\langle\bm{H},\bm{Y}_i\rangle]=&\wE[\wE[\pm\langle \bm{u},\bm{Y}_i\rangle]_{\bm{u}=\bm{H}}]\leq  \rho_i\wE[\Vert \bm{H}\Vert_2].
    \end{align}
    For non-negative $1$-dimensional function $G=G(\bm{Y_1,\cdots,\bm{Y}_{i-1}})\in\mathcal{H}$, we have 
    \begin{align}
    \label{eq:homo-estimate}
        \wE[G\Vert Y_i\Vert^2]=\wE[G]\wE[\Vert Y_i\Vert^2].
    \end{align}
    Now for fixed $\varepsilon>0$, denote 
    \begin{align*}
        &r_\varepsilon(\bm{x}):=(\Vert \bm{x}\Vert_2^2+\varepsilon^2)^{1/2}, \quad f_\varepsilon(\bm{x}):=r_\varepsilon^p(\bm{x}), \\
        &a_i:=\wE[\Vert \bm{Y}_i\Vert_2^p],\quad b_i:=\wE[\Vert \bm{Y}_i\Vert_2^2],\quad U_\varepsilon=\max_{0\le k\le n}\wE[r_\varepsilon^p(\bm{T}_k)].
    \end{align*}
    We first show that 
    \begin{align}
    \label{eq:Rosenthal-step1}
        U_\varepsilon\le C_p(\varepsilon^p+A+B^{p/2}+D^p).
    \end{align}
    Since $r_\varepsilon,f_\varepsilon\in C^{\infty}(\mathbb{R}^d)$, after some calculation, we can derive that 
    \begin{align*}
        & \nabla r_\varepsilon(\bm{x})=\frac{\bm{x}}{r_\varepsilon(\bm{x})}, \quad \nabla f_\varepsilon(\bm{x})=pr_\varepsilon^{p-2}(\bm{x})\bm{x},\\
        & D^2f_\varepsilon(\bm{x})=pr_\varepsilon^{p-2}(\bm{x})I+p(p-2)r_\varepsilon^{p-4}(\bm{x})\bm{x}\bm{x}^T.
    \end{align*}
    Note that the matrix $D^{2}f_{\varepsilon}(\bm{x})$ has an eigenvalue $p\,r_{\varepsilon}^{p-2}(\bm{x})$ of multiplicity $d-1$ and a simple eigenvalue $p\,r_{\varepsilon}^{p-4}(\bm{x})\bigl((p-1)\|\bm{x}\|_{2}^{2}+\varepsilon^{2}\bigr)$. The second derivative can be estimated as follows:
    \begin{align*}
        \Vert D^2f_\varepsilon(\bm{x})\Vert_{op}\leq \max\big\{pr_\varepsilon^{p-2}(\bm{x}), p\,r_{\varepsilon}^{p-4}(\bm{x})\bigl((p-1)\|\bm{x}\|_{2}^{2}+\varepsilon^{2}\bigr)\big\}\le p(p-1)r_\varepsilon^{p-2}(\bm{x}).
    \end{align*}
    Taylor’s formula with integral remainder gives
    \begin{align}
        f_\varepsilon(\bm{x}+\bm{y})&=f_\varepsilon (\bm{x}) +\langle \nabla f_\varepsilon(\bm{x}), \bm{y}\rangle+\int_0^1 (1-t)\langle \bm{y}, D^2f_\varepsilon(\bm{x}+t\bm{y})\bm{y}\rangle dt \nonumber\\
        &\le f_\varepsilon (\bm{x})+pr_\varepsilon^{p-2}(\bm{x})\langle\bm{x},\bm{y}\rangle+\int_0^1(1-t)\Vert D^2 f_\varepsilon(\bm{x}+t\bm{y})\Vert_{op}\Vert \bm{y}\Vert_2^2dt\nonumber\\
        & \leq f_\varepsilon (\bm{x})+pr_\varepsilon^{p-2}(\bm{x})\langle\bm{x},\bm{y}\rangle +\frac{p(p-1)}{2}\Vert \bm{y}\Vert_2^2\sup_{0\le t\le 1}r_\varepsilon^{p-2}(\bm{x}+t\bm{y})\nonumber\\
        &\le f_\varepsilon (\bm{x})+pr_\varepsilon^{p-2}(\bm{x})\langle\bm{x},\bm{y}\rangle  + C_p\Vert \bm{y}\Vert_2^2 \big(r_\varepsilon(\bm{x})+\Vert \bm{y}\Vert_2\big)^{p-2} \nonumber\\
        & \le f_\varepsilon (\bm{x})+pr_\varepsilon^{p-2}(\bm{x})\langle\bm{x},\bm{y}\rangle  + C_p\big(\Vert \bm{y}\Vert_2^2r_\varepsilon^{p-2}(\bm{x})+\Vert \bm{y}\Vert_2^p\big).
        \label{eq:tylor}
    \end{align}
    Substituting $\bm{x}=\bm{T}_{i-1}$ and $\bm{y}=\bm{Y}_i$ into \eqref{eq:tylor} and taking sublinear expectations gives
    \begin{equation}
    \label{eq:Tylar-expectation-1}
        \wE[f_\varepsilon(\bm{T}_i)]\le \wE[f_\varepsilon(\bm{T}_{i-1})]+p\wE[r_\varepsilon^{p-2}(\bm{T}_{i-1})\langle\bm{T}_{i-1},\bm{Y}_i\rangle]+C_p\wE[r_\varepsilon^{p-2}(\bm{T}_{i-1})\Vert \bm{Y}_i\Vert_2^2] +C_p\wE[\Vert \bm{Y}_i\Vert_2^p].
    \end{equation}
    Let $\bm{H}_{i-1}=r_\varepsilon^{p-2}(\bm{T}_{i-1}) \bm{T}_{i-1}$ and $G_{i-1}=r_\varepsilon^{p-2}(\bm{T}_{i-1})> 0$, by \eqref{eq:mean-estimate} and \eqref{eq:homo-estimate}, we have
    \begin{align}
        &\wE[\langle\bm{H}_{i-1},\bm{Y}_i\rangle]\leq \rho_i\wE[r_\varepsilon^{p-2}(\bm{T}_{i-1}) \Vert\bm{T}_{i-1}\Vert_2]\leq \rho_i \wE[r_\varepsilon^{p-1}(\bm{T}_{i-1}) ],\label{eq:H-estimate}\\
        &\wE[G_{i-1}\Vert \bm{Y}_i\Vert_2^2]=\wE[G_{i-1}]\wE[\Vert \bm{Y}_i\Vert_2^2]=b_i\wE[r_\varepsilon^{p-2}(\bm{T}_{i-1}) ].\label{eq:G-estimate}
    \end{align}
    Substituting \eqref{eq:H-estimate} and \eqref{eq:G-estimate} into \eqref{eq:Tylar-expectation-1}, and summing from $1$ to $k$ yields
    \begin{align*}
        \wE[f_\varepsilon(\bm{T}_k)]&\le \wE[f_\varepsilon(\bm{T}_{0})] +p\sum_{i=1}^k\rho_i \wE[r_\varepsilon^{p-1}(\bm{T}_{i-1}) ] +C_p\sum_{i=1}^k b_i\wE[r_\varepsilon^{p-2}(\bm{T}_{i-1}) ] +C_p\sum_{i=1}^ka_i\nonumber\\
        &\le \varepsilon^p + p\sum_{i=1}^k\rho_i \wE[r_\varepsilon^{p}(\bm{T}_{i-1}) ]^{1-1/p}+C_p\sum_{i=1}^k b_i\wE[r_\varepsilon^{p}(\bm{T}_{i-1}) ]^{1-2/p}+C_p\sum_{i=1}^ka_i.
    \end{align*}
    Taking maximum over $0\le k\le n$, it follows that
    \begin{equation*}
        U_\varepsilon\le \varepsilon^p+C_pA+C_pBU_\varepsilon^{1-2/p}+pDU_\varepsilon^{1-1/p}.
    \end{equation*}
    For $p=2$, the preceding estimate gives
$$
 U_\varepsilon
 \le \varepsilon^2+C(A+B)+2D U_\varepsilon^{1/2}
 \le \varepsilon^2+C(A+B)
      +\frac12 U_\varepsilon+2D^2.
$$
Absorbing $U_{\varepsilon}/2$ into the left-hand side gives \eqref{eq:Rosenthal-step1} when $p=2$.
For $p>2$, we apply Young's inequality as follows.
    Applying Lemma~\ref{lem:Young}, take $a=U_\varepsilon^{1-2/p}, b=B, r=p/(p-2),s=p/2$ and $\delta=1/(4C_p)$, we have
    \begin{equation*}
        C_pBU_\varepsilon^{1-2/p}\leq C_p(\delta U_\varepsilon + C_{r,s}B^{p/2})=\frac{1}{4}U_\varepsilon+C_pB^{p/2}.
    \end{equation*}
    Similarly, we can choose $\delta=1/(4p)$ such that $pDU_\varepsilon^{1-1/p}\le \frac14 U_\varepsilon+C_pD^p.$ Therefore, \eqref{eq:Rosenthal-step1} holds.
    Now we define 
    \begin{align*}
        & q:=p/(p-1),\quad r_k:=r_\varepsilon(\bm{T_k}),\quad  r_k^*=\max_{0\le j \le k}r_j,\\
        & \bm{H}_i:=(r_{i-1}^*)^{p-1}\frac{\bm{T}_{i-1}}{r_{i-1}},\quad Z_\varepsilon:=\wE[(r_n^*)^p].
    \end{align*}
    Lemma~\ref{lem:doob} gives that
    \begin{equation}
        (r_n^*)^p\leq q^pr_n^p-pq\sum_{i=1}^n(r_{i-1}^*)^{p-1}(r_i-r_{i-1}).
        \label{eq:max-estimate-2}
    \end{equation}
    Note that $r_\varepsilon(\bm{x})$ is convex, so $r_i-r_{i-1}\ge \langle \nabla r_\varepsilon(\bm{T}_{i-1}), \bm{Y}_i\rangle=\langle \frac{\bm{T}_{i-1}}{r_{i-1}},\bm{Y}_i\rangle$. Substituting it into \eqref{eq:max-estimate-2} and taking $\wE$ yields
    \begin{align*}
        Z_\varepsilon&\le q^p\wE[r_n^p]+pq\sum_{i=1}^n\wE[-\langle \bm{H}_i, \bm{Y}_i\rangle] \\
        & \le q^p U_\varepsilon+pq\sum_{i=1}^n\rho_i\wE[\Vert \bm{H}_i\Vert_2]\\
        & \le q^p U_\varepsilon+pq\sum_{i=1}^n\rho_i\wE[ (r_{i-1}^*)^{p}]^{1-1/p}\\
        & \le q^p U_\varepsilon+pqDZ_\varepsilon^{1-1/p} \le q^p U_\varepsilon+\frac12Z_\varepsilon+C_pD^p.
    \end{align*}
    The last inequality uses Lemma~\ref{lem:Young} again. Thus, combining \eqref{eq:Rosenthal-step1} obtains
    \begin{align*}
        \wE\big[\max_{0\le k\le n}\Vert \bm{T}_k\Vert_2^p\big]\le Z_\varepsilon\le C_p(U_\varepsilon+D^p)\le C_p\big(\varepsilon_p+A+B^{p/2}+D^p\big).
    \end{align*}
    Letting $\varepsilon\searrow0$ yields \eqref{eq:Rosenthal}.
    
    For general $\bm X_i$, we denote $\bm X_i^{(c)}=(X_{i,1}^{(c)},\ldots, X_{i,d}^{(c)})$, where $c>0$, $X^{(c)}=(-c)\vee X\wedge c$. With $\bm X_i^{(c)}$ in place of $\bm X_i$, define
$$
Y_{i,j}^{(c)}:=X_{i,j}^{(c)}-c_{i,j},\quad
\bm Y_i^{(c)}:=\bm X_i^{(c)}-\bm c_i,\quad
\bm T_0^{(c)}:=\bm 0,\quad
\bm T_k^{(c)}:=\sum_{i=1}^k\bm Y_i^{(c)},
\quad k=1,\cdots,n.
$$
Applying the bounded case of \eqref{eq:Rosenthal} to
$\bm X_1^{(c)},\ldots,\bm X_n^{(c)}$
yields the corresponding inequality for these truncated variables. Notice $\bm T_k=\bm T_k^{(c)}$, $k=1,\ldots, n$,  when   $|X_{i,j}|\le c$ for all $i=1,\ldots, n, j=1,\ldots, d$, and $ \max_{0\le k \le n}\Vert \bm{T}_k^{(c)}\Vert_2\le nd(c+\max_i\|\bm c_i\|_2)=:K_{c,n}$. Also,
    $$ \max_{0\le k \le n}\Vert \bm{T}_k\Vert_2^p\wedge K\le \max_{0\le k \le n}\Vert \bm{T}_k^{(c)}\Vert_2^p
   + K I\{ |X_{i,j}|>c \text{ for some } i, j\}, \;\; \forall K>0. $$
Since
$$ \V(|X_{i,j}|>c)\le \frac{\bE[|X_{i,j}|]}{c}\to 0 \text{ as } c\to \infty, $$
we have
$$\wE\bigg[\max_{0\le k \le n}\Vert \bm{T}_k\Vert_2^p\wedge K\bigg]\le \limsup_{c\to \infty}\wE\bigg[\max_{0\le k \le n}\Vert \bm{T}_k^{(c)}\Vert_2^p\bigg], \;\; \forall K>0. $$
Hence, 
$$\bE\bigg[\max_{0\le k \le n}\Vert \bm{T}_k\Vert_2^p\bigg]\le \limsup_{c\to \infty}\wE\bigg[\max_{0\le k \le n}\Vert \bm{T}_k^{(c)}\Vert_2^p\bigg]. $$
By \eqref{eq:Rosenthal}, it follows that 
$$\bE\bigg[\max_{0\le k \le n}\Vert \bm{T}_k\Vert_2^p\bigg]\le  C_p\limsup_{c\to\infty} \bigg\{\sum_{i=1}^n\wE[\|\bm Y_i^{(c)}\|_2^p]+\big(\sum_{i=1}^n\wE[\|\bm Y_i^{(c)}\|_2^2]\big)^{p/2}+\big(\sum_{i=1}^n\rho_i^c(\bm c_i)\big)^p\bigg\},$$
where $\rho_i^c(\bm c_i)=\max_{\bm \theta\in \Theta_i^c}\|\bm\theta-\bm c_i\|_2$, and $\Theta_i^c=\{\bm\theta\in\R^d:\langle\bm u, \bm\theta\rangle\le \wE[\langle \bm u, \bm X_i^{(c)}\rangle], \;\forall \bm u \in\R^d,\|\bm u\|_2\le 1\}$. 

Fix $i$ and let $c>\|\bm c_i\|_2$. Then $|c_{i,j}|<c$
for every $j=1,\ldots,d$.
If $|X_{i,j}|\le c$, then
$|Y_{i,j}^{(c)}|=|Y_{i,j}|$.
If $X_{i,j}>c$, then
$|Y_{i,j}^{(c)}|
=c-c_{i,j}
\le X_{i,j}-c_{i,j}
=|Y_{i,j}|.
$
Similarly, if $X_{i,j}<-c$, then
$
|Y_{i,j}^{(c)}|
=c+c_{i,j}
\le -X_{i,j}+c_{i,j}
=|Y_{i,j}|.
$
Thus, in all cases, $|Y_{i,j}^{(c)}|\le |Y_{i,j}|$, and hence
$\|\bm Y_i^{(c)}\|_2\le\|\bm Y_i\|_2$.
Moreover, the triangle inequality gives
$$
\|\bm Y_i^{(c)}\|_2
\le \|\bm X_i^{(c)}\|_2+\|\bm c_i\|_2
\le \sqrt{d}\,c+\|\bm c_i\|_2.
$$
Combining these bounds, we obtain
$$
\|\bm Y_i^{(c)}\|_2^p
\le
\|\bm Y_i\|_2^p
\wedge(\sqrt{d}\,c+\|\bm c_i\|_2)^p.
$$
Therefore, by monotonicity and the definition of $\bE$,
\begin{equation*}
\limsup_{c\to\infty}\wE[\|\bm Y_i^{(c)}\|_2^p]\le\lim_{c\to\infty}\wE\!\big[\|\bm Y_i\|_2^p\wedge(\sqrt{d}\,c+\|\bm c_i\|_2)^p\big] =\bE[\|\bm Y_i\|_2^p].
\end{equation*}
The same argument applies with $p$ replaced by $2$.

Also,
$$\left|\bE[\langle\bm{u},\bm{X}_i\rangle]-\wE[\langle\bm{u},\bm{X}_i^{(c)}\rangle]\right|
\le \|\bm u\|_2\bE[\|\bm X_i-\bm X_i^{(c)}\|_2]\le \|\bm u\|_2\frac{\bE[\|\bm X_i\|_2^2]}{c} \to 0 \text{ as } c\to \infty. $$
It follows that \eqref{eq:Rosenthal} holds for general $\bm X_i$s. 
    
\end{proof}

\begin{corollary}\label{cor:rosenthal}
    Suppose $\{\bm{X}_i;i\ge1\}$ and $\Theta_i,\rho_i$ are defined as in Theorem~\ref{thm:rosenthal}. Denote $\bm{X}_i=(X_{i,1},\cdots, X_{i,d})^\top$, and 
    \begin{align*}
        \overline{\bm\mu}_i:=\bE[\bm{X}_i]=\big(\bE[X_{i,1}],\cdots,\bE[X_{i,d}]\big)^\top,\quad \underline{\bm\mu}_i:=\be[\bm{X}_i]=\big(\be[X_{i,1}],\cdots,\be[X_{i,d}]\big)^\top.
    \end{align*}
    Then for any $i=1,\cdots,n$, we have 
    \begin{equation}\label{eq:rho-bound}
        \rho_i(\bm0)\leq \big(\Vert \overline{\bm\mu}_i\Vert_2^2+\Vert \underline{\bm\mu}_i\Vert_2^2\big)^{1/2}\leq \Vert \overline{\bm\mu}_i\Vert_2+\Vert \underline{\bm\mu}_i\Vert_2.
    \end{equation}
    In particular, for any $p\ge 2$, there exists a constant $C_p>0$ such that the following inequality holds:
    \begin{align}
        \bE\bigg[\max_{0\le k \le n}\Vert \sum_{i=1}^k\bm{X}_i\Vert_2^p\bigg]&\leq C_p\bigg\{\sum_{i=1}^n\bE[\Vert \bm{X}_i\Vert_2^p]+\bigg(\sum_{i=1}^n\bE[\Vert \bm{X}_i\Vert_2^2]\bigg)^{p/2}\nonumber\\
        &+\bigg[\sum_{i=1}^n\big(\Vert \bE[\bm{X}_i]\Vert_2+\Vert \be[\bm{X}_i]\Vert_2\big)\bigg]^p\bigg\}.
    \end{align}
\end{corollary}

\begin{proof}
Let $\bm{e}_1,\dots,\bm{e}_d$ denote the standard orthonormal basis of $\R^{d}$. By the definition of
$\Theta_i$, for every $\bm\theta\in\Theta_i$ and every $j=1,\dots,d$ the choices $\bm{u}=\bm{e}_j$ and
$\bm{u}=-\bm{e}_j$ give
\begin{align}
\theta_j=\langle \bm{e}_j,\bm\theta\rangle
&\le\bE\langle \bm{e}_j,\bm{X}_i\rangle
=\bE\bigl[X_{i,j}\bigr]
=\overline\mu_{i,j}, \label{eq:upper-theta}\\
-\theta_j=\langle -\bm{e}_j,\bm\theta\rangle
&\le\bE\langle -\bm{e}_j,\bm{X}_i\rangle
=\bE\bigl[-X_{i,j}\bigr]
=-\be\bigl[X_{i,j}\bigr]
=-\underline\mu_{i,j} .\label{eq:lower-theta}
\end{align}
Combining \eqref{eq:upper-theta} and \eqref{eq:lower-theta} yields
\begin{equation}
\underline\mu_{i,j}\le\theta_j\le\overline\mu_{i,j},
\qquad\text{hence}\qquad
|\theta_j|^{2}\le\max\bigl\{|\underline\mu_{i,j}|^{2},|\overline\mu_{i,j}|^{2}\bigr\}
\label{eq:theta-bound}
\end{equation}
for every $j=1,\dots,d$ and every $\theta\in\Theta_i$.

Since the right-hand side of \eqref{eq:theta-bound} does not depend on $\theta$, taking the supremum
over $\theta\in\Theta_i$ gives
\begin{align}
\rho_i(\bm{0})^{2}
&=\max_{\bm\theta\in\Theta_i}\|\bm\theta\|_2^{2}
=\max_{\bm\theta\in\Theta_i}\sum_{j=1}^{d}|\theta_j|^{2}\le\sum_{j=1}^{d}\max\bigl\{|\underline\mu_{i,j}|^{2},|\overline\mu_{i,j}|^{2}\bigr\}\nonumber\\
&\le\sum_{j=1}^{d}\Bigl(|\underline\mu_{i,j}|^{2}+|\overline\mu_{i,j}|^{2}\Bigr)
=\|\underline{\bm\mu}_i\|_2^{2}+\|\overline{\bm\mu}_i\|_2^{2}\le\bigl(\|\underline{\bm\mu}_i\|_2+\|\overline{\bm\mu}_i\|_2\bigr)^{2}.
\label{eq:rho2}
\end{align}
Taking square roots in \eqref{eq:rho2} proves \eqref{eq:rho-bound}.

Finally, applying Theorem~\ref{thm:rosenthal} with
the centering vectors $\bm{c}_i=\bm0$, and
$\rho_i(\bm0)=\max_{\bm\theta\in\Theta_i}\|\bm\theta\|_2$, we obtain
\begin{align*}
\bE\Bigl[\max_{0\le k\le n}\Bigl\|\sum_{i=1}^{k}\bm{X}_i\Bigr\|_2^{p}\Bigr]
&\le C_p\Bigl\{\sum_{i=1}^{n}\bE\|\bm{X}_i\|_2^{p}
+\Bigl(\sum_{i=1}^{n}\bE\|\bm{X}_i\|_2^{2}\Bigr)^{p/2}
+\Bigl(\sum_{i=1}^{n}\rho_i(\bm0)\Bigr)^{p}\Bigr\}\\
&\le C_p\Bigl\{\sum_{i=1}^{n}\bE\|\bm{X}_i\|_2^{p}
+\Bigl(\sum_{i=1}^{n}\bE\|\bm{X}_i\|_2^{2}\Bigr)^{p/2}
+\Bigl[\sum_{i=1}^{n}\bigl(\|\bE[\bm{X}_i]\|_2+\|\be[\bm{X}_i]\|_2\bigr)\Bigr]^{p}\Bigr\}.
\end{align*}
\end{proof}

\section{Well-posedness under weak moment assumptions}
\label{sec:weak-wellposedness}

Throughout this section, we show the well-posedness of the linear process on the canonical space under the following consequences of Assumption~\ref{ass:innovations}:
\begin{equation}\label{eq:wp-weak-moments}
 \bE[\varepsilon_1]
 =
 \be[\varepsilon_1]
 =0,
 \qquad
 \bar\sigma^2
 := \bE[\varepsilon_1^2] =
 \sup_{c>0}
 \widehat{\mathbb E}[\varepsilon_1^2\wedge c]
 <\infty.
\end{equation}
The condition
$x^2\mathbb V(|\varepsilon_1|\ge x)\to0$
is not needed for the construction below.

We first construct a copy of the innovation sequence and introduce the associated moment spaces.
\begin{lemma}\label{lem:wp-canonical-representation}
Let
\[
 \Omega^{\mathrm c}:=\mathbb R^{\mathbb Z},
 \qquad
 \mathcal F^{\mathrm c}
 :=\mathcal B(\mathbb R)^{\otimes\mathbb Z},
 \qquad
 \xi_j(\omega):=\omega_j,\quad j\in\mathbb Z.
\]
There exists a nonempty weakly compact family $\mathcal P$
of probability measures $P$ on
$(\Omega^{\mathrm c},\mathcal F^{\mathrm c})$ such that
\begin{equation}\label{eq:wp-cylinder-expectation}
 \sup_{P\in\mathcal P}E_P
 [\varphi(\xi_{j_1},\ldots,\xi_{j_m})]
 =
 \widehat{\mathbb E}
 [\varphi(\varepsilon_{j_1},\ldots,\varepsilon_{j_m})]
\end{equation}
for every $j_1<\cdots<j_m$ and
$\varphi\in C_{b,\mathrm{Lip}}(\mathbb R^m)$.
Moreover,
\begin{equation}\label{eq:wp-canonical-second-moment}
 \sup_{P\in\mathcal P}E_P[\xi_j^2]
 =\bar\sigma^2,
 \qquad j\in\mathbb Z.
\end{equation}
The set function
\begin{equation}\label{eq:wp-capacity}
 \mathbb V^{\mathcal P}(A)
 :=\sup_{P\in\mathcal P}P(A),
 \qquad A\in\mathcal F^{\mathrm c},
\end{equation}
is a countably subadditive capacity. 
\end{lemma}

\begin{proof}
For every $j\in\mathbb Z$ and $M>0$,
\[
 \mathbb V(|\varepsilon_j|\ge M)
 \le M^{-2}\widehat{\mathbb E}
       [\varepsilon_j^2\wedge M^2]
 \le\bar\sigma^2M^{-2}
 \longrightarrow0.
\]
Thus each $\varepsilon_j$ is tight. Choose a bijection
$\iota:\mathbb N\to\mathbb Z$ and apply
\cite[Proposition~2.1(a)--(c)]{Zhang2023SLLN}
to $(\varepsilon_{\iota(k)})_{k\ge1}$.
Identifying the resulting product space
$\mathbb R^{\infty}$ with $\mathbb R^{\mathbb Z}$
via $\iota$ gives \eqref{eq:wp-cylinder-expectation}. Nonemptiness follows by applying the representation
in Proposition~2.1(c) to the constant function $1$.

By monotone convergence,
\begin{equation*}
    \sup_{P\in\mathcal P}E_P[\xi_j^2]=\sup_{P\in\mathcal P}\sup_{a>0}
   E_P[\xi_j^2\wedge a]=\sup_{a>0}\sup_{P\in\mathcal P}
   E_P[\xi_j^2\wedge a]=\sup_{a>0}\widehat{\mathbb E}
   [\varepsilon_j^2\wedge a]
 =\bar\sigma^2.
\end{equation*}
This proves \eqref{eq:wp-canonical-second-moment}.

The assertions concerning $\mathbb V^{\mathcal P}$ follow
directly from
\cite[Proposition~2.1(b)--(c)]{Zhang2023SLLN}.
\end{proof}

Define
\[
 \mathcal H^{\mathrm c}
 :=
 \left\{
 Z:\Omega^{\mathrm c}\to\mathbb R:
 Z \text{ is }\mathcal F^{\mathrm c}\text{-measurable},
 \quad
 \sup_{P\in\mathcal P}E_P[|Z|]<\infty
 \right\},
\]
and equip this space with the upper expectation
\begin{equation}\label{eq:wp-representation}
 \mathbb E^{\mathrm c}[Z]
 :=
 \sup_{P\in\mathcal P}E_P[Z],
 \qquad Z\in\mathcal H^{\mathrm c}.
\end{equation}
The space $\mathcal H^{\mathrm c}$ is a vector lattice
containing the constants, and $\mathbb E^{\mathrm c}$
is finite-valued, monotone, constant-preserving,
subadditive, and positively homogeneous.
Moreover, $\mathcal H^{\mathrm c}$ contains all bounded
measurable functions and is therefore stable under
bounded Lipschitz compositions. Thus
$(\Omega^{\mathrm c},\mathcal H^{\mathrm c},
\mathbb E^{\mathrm c})$
is a sublinear expectation space; see
\cite[Sections~1.4 and~6.1.2]{Peng2019}.

For $p\ge1$ and a measurable random vector
$\bm U:\Omega^{\mathrm c}\to\mathbb R^d$, write
\[
 \|\bm U\|_{p,\mathcal P}
 :=
 \left(
 \sup_{P\in\mathcal P}E_P[\|\bm U\|_2^p]
 \right)^{1/p}.
\]
Let $L_{\mathcal P}^p(\mathbb R^d)$ be the space of
such vectors with finite $\|\cdot\|_{p,\mathcal P}$,
modulo quasi-sure equality.
For scalar variables, we abbreviate
$L_{\mathcal P}^p(\mathbb R)$ to $L_{\mathcal P}^p$.

We also write: 
\begin{equation}\label{eq:wp-canonical-expectation}
 \mathcal E^{\mathrm c}[Z]
 :=
 -\mathbb E^{\mathrm c}[-Z],
 \qquad Z\in L_{\mathcal P}^1.
\end{equation}
By \eqref{eq:wp-canonical-second-moment},
each $\xi_j$ belongs to $L_{\mathcal P}^2$.

The following conditional moment bounds yield a second-moment estimate for finite weighted sums. Let
\[
 \mathcal F_j^{\mathrm c}
 :=\sigma(\xi_k:k\le j),\qquad j\in\mathbb Z.
\]

\begin{lemma}\label{lem:wp-conditional-moments}
For every $P\in\mathcal P$ and $j\in\mathbb Z$,
\begin{equation}\label{eq:wp-martingale-moments}
 E_P[\xi_j\mid\mathcal F_{j-1}^{\mathrm c}]=0,
 \qquad
 E_P[\xi_j^2\mid\mathcal F_{j-1}^{\mathrm c}]
 \le\bar\sigma^2,
 \quad P\text{-a.s.}
\end{equation}
Consequently, for every finite $I\subset\mathbb Z$
and deterministic $\bm{v}_j\in\mathbb R^d$,
\begin{equation}\label{eq:wp-finite-L2}
 \mathbb E^{\mathrm c}
 \left[
 \left\|\sum_{j\in I}\bm{v}_j\xi_j\right\|_2^2
 \right]
 \le
 \bar\sigma^2\sum_{j\in I}\|\bm{v}_j\|_2^2.
\end{equation}
\end{lemma}

\begin{proof}
Fix $P\in\mathcal P$ and $j\in\mathbb Z$. Put
\[
 \mathcal G_{j,m}:=\sigma(\xi_{j-m},\ldots,\xi_{j-1}),
 \qquad
 \mathcal G_{j,m}\uparrow\mathcal F_{j-1}^{\mathrm c}.
\]
By \eqref{eq:wp-cylinder-expectation} and
\eqref{eq:wp-representation}, the sequence
$\xi_{j-m},\ldots,\xi_j$ is independent
 under $\mathbb E^{\mathrm c}$ for bounded
Lipschitz test functions.
Applying Proposition~2.4 of 
\cite{GuoLiLi2022MeanUncertainty} therefore
gives, for every $\varphi\in C_{b,\mathrm{Lip}}(\mathbb R)$,
\[
 E_P[\varphi(\xi_j)\mid\mathcal G_{j,m}]
 \le\mathbb E^{\mathrm c}[\varphi(\xi_j)]
 =\widehat{\mathbb E}[\varphi(\varepsilon_1)],
 \qquad P\text{-a.s.}
\]
Letting $m\to\infty$ by the upward martingale
convergence theorem gives
\begin{equation*}\label{eq:wp-conditional-domination}
 E_P[\varphi(\xi_j)\mid\mathcal F_{j-1}^{\mathrm c}]
 \le\widehat{\mathbb E}[\varphi(\varepsilon_1)],
 \qquad P\text{-a.s.}
\end{equation*}

Taking $\varphi(x)=\pm x^{(c)}$ and $\varphi(x)=x^2\wedge b$,
respectively, yields
\[
\begin{aligned}
 &-\widehat{\mathbb E}[-\varepsilon_1^{(c)}]
 \le E_P[\xi_j^{(c)}\mid\mathcal F_{j-1}^{\mathrm c}]
 \le\widehat{\mathbb E}[\varepsilon_1^{(c)}],\\
 &E_P[\xi_j^2\wedge b\mid\mathcal F_{j-1}^{\mathrm c}]
 \le\widehat{\mathbb E}[\varepsilon_1^2\wedge b]
 \le\bar\sigma^2.
\end{aligned}
\]
Moreover, by \eqref{eq:wp-canonical-second-moment},
\[
 E_P[|\xi_j-\xi_j^{(c)}|]
 \le c^{-1}E_P[\xi_j^2]
 \le\frac{\bar\sigma^2}{c}\longrightarrow0.
\]
Letting $c\to\infty$ and 
$b\to\infty$ yields
\eqref{eq:wp-martingale-moments}.

For $i<j$, square integrability and
\eqref{eq:wp-martingale-moments} give
\[
 E_P[\xi_i\xi_j] = E_P[E_P[\xi_i \xi_j\mid \mathcal F_{j-1}^{\mathrm c}]]
 =E_P\left[\xi_i
   E_P[\xi_j\mid\mathcal F_{j-1}^{\mathrm c}]\right]
 =0.
\]
Consequently,
\[
\begin{aligned}
 E_P\left[\left\|\sum_{j\in I}\bm{v}_j\xi_j\right\|_2^2\right]
 =\sum_{j\in I}\|\bm{v}_j\|_2^2E_P[\xi_j^2]
 \le\bar\sigma^2\sum_{j\in I}\|\bm{v}_j\|_2^2.
\end{aligned}
\]
Taking the supremum over $P\in\mathcal P$ proves
\eqref{eq:wp-finite-L2}.
\end{proof}

This estimate allows us to pass from finite sums to series with square-summable coefficients.

\begin{theorem}\label{thm:wp-weighted-series}
Let $\bm v=(\bm v_j)_{j\in\mathbb Z}
\in\ell^2(\mathbb Z;\mathbb R^d)$, that is,
$\bm v_j\in\mathbb R^d$ and $
 \|\bm v\|_{\ell^2}^2
 :=\sum_{j\in\mathbb Z}\|\bm v_j\|_2^2<\infty.$
Then there exists a measurable vector
$\bm{\mathcal I}(\bm{v})\in L_{\mathcal P}^2(\mathbb R^d)$ such that
\[
 \sum_{|j|\le N}\bm{v}_j\xi_j
 \longrightarrow\bm{\mathcal I}(\bm{v}) 
 \quad\text{quasi-surely and in }L_{\mathcal P}^2(\mathbb R^d).
\]
The limit is unique up to quasi-sure equality and satisfies
\begin{equation}\label{eq:wp-series-tail}
 \mathbb E^{\mathrm c}
 \left[
 \left\|
 \bm{\mathcal I}(\bm{v})-\sum_{j\in I}\bm{v}_j\xi_j
 \right\|_2^2
 \right]
 \le
 \bar\sigma^2\sum_{j\notin I}\|\bm{v}_j\|_2^2
\end{equation}
for every finite $I\subset\mathbb Z$. In particular, 
\begin{equation}\label{eq:limit-l2-norm}
 \left\|
 \bm{\mathcal I}(\bm{v})
 \right\|_{2,\mathcal{P}}^2
 \le
 \bar\sigma^2\|\bm v\|_{\ell^2}^2.
\end{equation}

The mapping $\bm{\mathcal I}:\ell^2(\mathbb Z;\mathbb R^d)\longrightarrow L_{\mathcal P}^2(\mathbb R^d)$ is linear.
Also, for each $P\in\mathcal P$,
\[
 E_P[\bm{\mathcal I}(\bm{v})]=\bm{0}.
\]
\end{theorem}

\begin{proof}
Write
\[
 \bm{T}_N:=\sum_{|j|\le N}\bm{v}_j\xi_j.
\]
Fix $M>N$. Set $\bm U_{-M-1}:=\bm0$ and define
\[
\bm U_k:=\sum_{j=-M}^{k}
\mathbf1_{\{|j|>N\}}\bm v_j\xi_j,
\qquad -M\le k\le M.
\]
For $N<k\le M$, we have
$\bm T_k-\bm T_N=\bm U_k-\bm U_{-k-1}$.
By Corollary \ref{cor:rosenthal}, it follows that, 
\begin{align}\label{eq:wp-cauchy-L2}
&\E^c [\max_{N< k\le M}\|\bm{T}_k-\bm{T}_N\|_2^2]\le 4\E^c[\max_{-M-1\le k\le M}\|\bm U_k\|_2^2]\nonumber\\
=&4\sup_{P\in\mathcal{P}}E_P\bigg[\max_{-M-1< k\le M}\bigg\|\sum_{-M-1<j\le k}\mathbf1_{\{|j|> N\}}\bm{v}_j\xi_j\bigg\|_2^2\bigg]
 \le C\bar\sigma^2\sum_{N<|j|\le M}\|\bm{v}_j\|_2^2\longrightarrow0
\end{align}
as $M>N\to\infty$. 
Choose $N_k\uparrow\infty$ such that
\[
 \sup_{P\in\mathcal{P}}E_P\bigg[\max_{N_k< N\le N_{k+1}}\bigg\|\sum_{N_k< |j|\le N}\bm{v}_j\xi_j\bigg\|_2^2\bigg]
 \le2^{-4k}.
\]
Then
\begin{align*}
 \sum_{k=1}^\infty
 \mathbb V^{\mathcal P}\left(
 \max_{N_k< N\le N_{k+1}}\|\bm{T}_{N}-\bm{T}_{N_k}\|_2>2^{-k}
 \right)
 \le
 \sum_{k=1}^\infty
 2^{2k}
 \E^c\left[\max_{N_k<N\le N_{k+1}}\|\bm{T}_{N}-\bm{T}_{N_k}\|^2_2\right]\le\sum_{k=1}^\infty2^{-2k}<\infty.
\end{align*}
Since $\mathbb V^{\mathcal P}$ is countably subadditive,
\[
 \mathbb V^{\mathcal P}\left(
 \limsup_{k\to\infty}
 \{\max_{N_k<N\le N_{k+1}}\|\bm{T}_{N}-\bm{T}_{N_k}\|^2>2^{-k}\}
 \right)=0.
\]
Thus $\bm{T}_N$ converges outside a polar set.
Define $\bm{\mathcal I}(\bm{v})$ to be this limit there, and
set it equal to $\bm{0}$ on the exceptional set.
It is an $\mathcal F^{\mathrm c}$-measurable random vector.

For every $P\in\mathcal P$, since $T_N\to \mathcal{I}(\bm{v}) \enspace a.s. P$. Fatou's lemma gives
\begin{align*}
 E_P[\|\bm{T}_N-\bm{\mathcal I}(\bm{v})\|_2^2] &= E_P[\liminf_{M\to\infty}\|\bm{T}_N-\bm{T}_{M}\|_2^2]
 \le
 \liminf_{M\to\infty}
 E_P[\|\bm{T}_N-\bm{T}_{M}\|_2^2]\\
 & \le \liminf_{M\to\infty} \|\bm{T}_N-\bm{T}_{M}\|_{2,\mathcal{P}}^2\le
 \bar\sigma^2
 \sum_{|j|>N}\|\bm{v}_j\|_2^2.
\end{align*}
Taking the supremum over $P$ yields
\[
 \|\bm{T}_N-\bm{\mathcal I}(\bm{v})\|_{2,\mathcal P}^2
 \le
 \bar\sigma^2
 \sum_{|j|>N}\|\bm{v}_j\|_2^2
 \longrightarrow0.
\]
Similarly, for every finite $I\subset\mathbb Z$,
\begin{align*}
 E_P\left[
 \left\|\bm{\mathcal I}(\bm{v})-\sum_{j\in  I}\bm{v}_j\xi_j\right\|_2^2
 \right]
 \le
 \liminf_{k\to\infty}
 E_P\left[
 \left\|\bm{T}_{N_k}-\sum_{j\in  I}\bm{v}_j\xi_j\right\|_2^2
 \right]\le
 \bar\sigma^2\sum_{j\notin  I}\|\bm{v}_j\|_2^2.
\end{align*}
This proves \eqref{eq:wp-series-tail}, and setting $I=\varnothing$ yields \eqref{eq:limit-l2-norm}.

Uniqueness follows from
\[
 \|\bm{U}-\bm{V}\|_{2,\mathcal P}=0
 \Longrightarrow
 \mathbb V^{\mathcal P}(\|\bm{U}-\bm{V}\|_2>\delta)=0
 \quad \text{for any }\delta>0
 \Longrightarrow \bm{U}=\bm{V}\quad\text{quasi-surely}.
\]
Linearity follows by passing to the limits of finite sums.
Finally, since $
 E_P[\bm{T}_N]=\bm{0},
 E_P[\|\bm{\mathcal I}(\bm{v})-\bm{T}_N\|_2]
 \le\|\bm{\mathcal I}(\bm{v})-\bm{T}_N\|_{2,\mathcal P}\longrightarrow0,
$ it obtains that $E_P[\bm{\mathcal I}(\bm{v})]=\bm{0}$.
\end{proof}

For scalar coefficient sequences, we write $\mathcal I$ for
the scalar-valued specialization of $\bm{\mathcal I}$.

\begin{corollary}\label{cor:wp-linear-process}
Suppose that $\bm{\alpha}=(\alpha_j)_{j\in\mathbb Z}\in\ell^2(\mathbb Z)$.
Then
\[
 Y_i^{\mathrm c}
 :=\sum_{j\in\mathbb Z}\alpha_j\xi_{i-j},
 \qquad i\in\mathbb Z,
\]
is well defined quasi-surely and in $L_{\mathcal P}^2$, and
\begin{equation}\label{eq:wp-linear-moments}
 \mathbb E^{\mathrm c}[(Y_i^{\mathrm c})^2]
 \le\bar\sigma^2\sum_{j\in\mathbb Z}\alpha_j^2,
 \qquad
 \mathbb E^{\mathrm c}[Y_i^{\mathrm c}]
 = \mathcal{E}^{\mathrm c}[Y_i^{\mathrm c}]
 =0.
\end{equation}
For $m\ge1$,
\begin{equation}\label{eq:wp-partial-sums}
 S_m^{\mathrm c}
 :=\sum_{i=1}^mY_i^{\mathrm c}
 =
 \sum_{j\in\mathbb Z}
 \left(\sum_{i=1}^m\alpha_{i-j}\right)\xi_j
 \quad\text{ quasi-surely and in }L_{\mathcal P}^2.
\end{equation}
\end{corollary}

\begin{proof}
For $\bm{a}_{i}:=(\alpha_{i-j})_{j\in\mathbb Z},$
we have
$$
 \|\bm{a}_{i}\|_{\ell^2}^2
 =\sum_{j\in\mathbb Z}\alpha_j^2<\infty.
$$
Hence $Y_i^{\mathrm c}=\mathcal I(\bm{a}_{i})$ exists,
and \eqref{eq:wp-linear-moments} follows from
Theorem~\ref{thm:wp-weighted-series}.

Moreover,
$$
 \sum_{j\in\mathbb Z}
 \left|\sum_{i=1}^m\alpha_{i-j}\right|^2
 \le
 m\sum_{i=1}^m\sum_{j\in\mathbb Z}\alpha_{i-j}^2
 =m^2\|\bm{\alpha}\|_{\ell^2}^2<\infty.
$$
Hence $\sum_{i=1}^{m}\bm{a}_i\in\ell^2(\mathbb{Z})$. By linearity,
\[
 S_m^{\mathrm c}
 =\sum_{i=1}^m\mathcal I(\bm{a}_{i})
 =\mathcal I\left(\sum_{i=1}^m \bm{a}_{i}\right),
\]
which proves \eqref{eq:wp-partial-sums}.
\end{proof}

We finally relate the copied process to finite truncations on the original sublinear expectation space through bounded Lipschitz test functions. For $N\ge1$, define
\[
 Y_i^{(N)}
 :=\sum_{|j|\le N}\alpha_j\varepsilon_{i-j},
 \qquad
 Y_i^{\mathrm c,(N)}
 :=\sum_{|j|\le N}\alpha_j\xi_{i-j}.
\]
By \eqref{eq:wp-cylinder-expectation} and
\eqref{eq:wp-representation}, for every
$\varphi\in C_{b,\mathrm{Lip}}(\mathbb R^d)$,
\begin{align*}
 \widehat{\mathbb E}
 \bigl[\varphi(Y_{i_1}^{(N)},\ldots,Y_{i_d}^{(N)})\bigr]=
 \mathbb E^{\mathrm c}
 \bigl[\varphi(
 Y_{i_1}^{\mathrm c,(N)},\ldots,
 Y_{i_d}^{\mathrm c,(N)})\bigr].
\end{align*}
Furthermore,
\begin{equation*}
\mathbb E^c\!\left[
 \sum_{r=1}^{d}|Y_{i_r}^{c,(N)}-Y_{i_r}^{c}|^2
\right]\le
\sum_{r=1}^{d}\mathbb E^c\!\left[
 \left|\sum_{|j|>N}\alpha_j\xi_{i_r-j}\right|^2
\right]\le d\bar\sigma^2\sum_{|j|>N}\alpha_j^2.
\end{equation*}
Consequently,
\begin{equation}\label{eq:wp-transfer-fdd}
\begin{aligned}
 &\left|
 \widehat{\mathbb E}
 \bigl[\varphi(Y_{i_1}^{(N)},\ldots,Y_{i_d}^{(N)})\bigr]
 -
 \mathbb E^{\mathrm c}
 \bigl[\varphi(Y_{i_1}^{\mathrm c},\ldots,
              Y_{i_d}^{\mathrm c})\bigr]
 \right|\le 
 C\bar\sigma
 \left(d\sum_{|j|>N}\alpha_j^2\right)^{1/2}
 \longrightarrow0.
\end{aligned}
\end{equation}

For the linear interpolation, put
\[
 \widetilde X_n^{\mathrm c}(t)
 :=
 n^{-H}\left[
 S_{\lfloor nt\rfloor}^{\mathrm c}
 +(nt-\lfloor nt\rfloor)
 Y_{\lfloor nt\rfloor+1}^{\mathrm c}
 \right].
\]
For each fixed $n$, this is a measurable
$C[0,1]$-valued stochastic process. Also,
\[
 \|\widetilde X_n^{\mathrm c}\|_\infty
 \le n^{-H}\sum_{i=1}^n|Y_i^{\mathrm c}|,
\]
and hence
\begin{equation}\label{eq:wp-polygonal-moment}
 \mathbb E^{\mathrm c}
 [\|\widetilde X_n^{\mathrm c}\|_\infty^2]
 \le n^{1-2H}\sum_{i=1}^n\mathbb{E}^{\mathrm{c}}[(Y_i^c)^2]\le 
 \bar\sigma^2n^{2-2H}\|\bm{\alpha}\|_{\ell^2}^2
 <\infty.
\end{equation}

Let $\widetilde X_{n,N}$ and
$\widetilde X_{n,N}^{\mathrm c}$ be the corresponding
interpolations obtained from
$Y_i^{(N)}$ and $Y_i^{\mathrm c,(N)}$, respectively.
Then
\[
 \mathbb E^{\mathrm c}
 [\|\widetilde X_{n,N}^{\mathrm c}
       -\widetilde X_n^{\mathrm c}\|_\infty^2]
 \le
 \bar\sigma^2n^{2-2H}
 \sum_{|j|>N}\alpha_j^2 \longrightarrow0\qquad\text{as } N\to\infty.
\]
Therefore, for fixed $n$ and
$\Phi\in C_{b,\mathrm{Lip}}(C[0,1])$,
\begin{equation}\label{eq:wp-transfer-path}
 \lim_{N\to\infty}
 \widehat{\mathbb E}[\Phi(\widetilde X_{n,N})]
 = \lim_{N\to\infty}\mathbb E^{\mathrm c}[\Phi(\widetilde X^{\mathrm c}_{n,N})]=
 \mathbb E^{\mathrm c}[\Phi(\widetilde X_n^{\mathrm c})].
\end{equation}

If the extension of $\wE$ on bounded Lipschitz functions of $Y_1,\ldots, Y_n$ is defined by \eqref{eq:extestionE:1}, then for any $\varphi\in C_{b,Lip}(\mathbb R^n)$ and $\Phi\in C_{b,Lip}(C[0,1])$, we have
$$ \wE[\varphi(Y_1,\ldots,Y_n)]=\E^c[\varphi(Y_1^c,\ldots,Y_n^c)]
=\sup_{P\in \mathcal{P}}E_P[\varphi(Y_1^c,\ldots,Y_n^c)], $$
$$\wE[\Phi(\widetilde X_n)]=\mathbb E^{\mathrm c}[\Phi(\widetilde X_n^{\mathrm c})]
=\sup_{P\in \mathcal{P}}E_P[\Phi(\widetilde X_n^{\mathrm c})]. $$

Now, for any $\varphi\in C_b(\mathbb R^n)$ and a sequence $\varphi_k\in C_{b, Lip}(\mathbb R^n)$ with
$$\|\varphi-\varphi_k\|_{[-R,R]}:=\sup_{x_1,\ldots, x_n\in[-R, R]}|\varphi(x_1,\ldots,x_n)-\varphi_k(x_1,\ldots, x_n)|
\to 0 \;\text{ as } k\to \infty,\; \forall R>0, $$
without loss of generality, we can assume $|\varphi(x_1,\ldots,x_n)|\le 1$, $|\varphi_k(x_1,\ldots,x_n)|\le 1$. Thus,
$$ \E^c|\varphi(Y_1^c,\ldots,Y_n^c)-\varphi_k(Y_1^c,\ldots, Y^c_n)|\le 
\|\varphi-\varphi_k\|_{[-R.R]}+2\mathbb V^{\mathcal{P}}\Big(\bigcup_{i=1}^n |Y_i^c|\ge R\Big)\to 0$$
as $k\to \infty$ and then $R\to \infty$.  Thus, 
$$\wE[\varphi_k(Y_1,\ldots, Y_n)]=\E^c[\varphi_k(Y_1^c,\ldots, Y^c_n)]\to \E^c[\varphi(Y_1^c,\ldots, Y^c_n)]. $$
It follows that the extension of $\wE$ to that on bounded continuous functions of $Y_1,\ldots, Y_n$ is well-defined by \eqref{eq:extestionE:2}, and
$$ \wE[\varphi(Y_1,\ldots,Y_n)]=\E^c[\varphi(Y_1^c,\ldots,Y_n^c)]
=\sup_{P\in \mathcal{P}}E_P[\varphi(Y_1^c,\ldots,Y_n^c)], \;\; \varphi\in C_b(\mathbb R^n), $$
$$\wE[\Phi(\widetilde X_n)]=\mathbb E^{\mathrm c}[\Phi(\widetilde X_n^{\mathrm c})]
=\sup_{P\in \mathcal{P}}E_P[\Phi(\widetilde X_n^{\mathrm c})],  \;\; \Phi\in C_b(C[0,1]).$$

\begin{remark}\label{rem:wp-scope}
The construction preserves the bounded Lipschitz
finite-dimensional distributions of the innovations
and their truncated variance parameters.
In the sequel, we work with the copied process in
$L_{\mathcal P}^2$ under $\mathbb E^{\mathrm c}$;
its relation to the original finite truncations is
given by \eqref{eq:wp-transfer-fdd} and
\eqref{eq:wp-transfer-path}, or equivalently, \eqref{eq:extestionE:1} and \eqref{eq:extestionE:2}. 
  From now on, we write
$\varepsilon_j$, $Y_i$, $S_m$, and $\widehat{\mathbb E}$
for the copied quantities
$\xi_j$, $Y_i^{\mathrm c}$, $S_m^{\mathrm c}$,
and $\mathcal E^{\mathrm c}$, respectively.
\end{remark}

\begin{remark}
\label{rem:wp-CC}
If $\wE$ satisfies condition (CC), let $\mathcal P_0$ be a family
of probability measures on $(\Omega,\sigma(\mathcal H))$
representing $\widehat{\mathbb E}$ on its bounded
variables.
Set
$\mathbb E_{\mathcal P_0}[Z]:=\sup_{P\in\mathcal P_0}E_P[Z]$
on $L^1_{\mathcal P_0}$, with the moment spaces
defined as above.
The preceding arguments apply on this space and give
\[
 Y_i^0
 :=\lim_{N\to\infty}
 \sum_{|j|\le N}\alpha_j\varepsilon_{i-j}
 \quad\text{in }L^2_{\mathcal P_0}.
\]
The resulting process has the same bounded Lipschitz
finite-dimensional distributions as the copied process.
Truncating the innovations also shows that $Y_i^0$
belongs to the $L^1_{\mathcal P_0}$-closure of
$\mathcal H_b$.
On this closure, $\mathbb E_{\mathcal P_0}$ is the
continuous extension of
$\widehat{\mathbb E}|_{\mathcal H_b}$.
\end{remark}

\section{Finite-dimensional convergence}\label{sec:fdc}
We first show that the rescaled partial-sum coefficients
approximate the limiting kernel uniformly on bounded intervals. Write
\begin{equation*}
b_{n,t,j}:=\sum_{i=1}^{\lfloor nt\rfloor}\alpha_{i-j}.
\end{equation*}

\begin{lemma}\label{lem:local-profile}
Under Assumption~\ref{ass:coeff}, for every fixed $t\in[0,1]$ and $R>0$,
\begin{equation}
\sup_{|j|\le Rn}\bigl|n^{\beta-1}b_{n,t,j}-f_t(j/n)\bigr|\longrightarrow0 \quad \text{as }n\to\infty  .
\label{eq:coefapprox}
\end{equation}
\end{lemma}

\begin{proof}
The assertion is immediate for $t=0$, so suppose $t>0$. 
Assumption~\ref{ass:coeff}
implies the global bound
\begin{equation}
|\alpha_r|\le C(1+|r|)^{-\beta},\qquad r\in\Z.
\label{eq:alpha-global-bound}
\end{equation}  
Fix $\delta>0$ and $j$, decompose $b_{n,t,j}=b_{n,t,j}^{1,\delta}+b_{n,t,j}^{2,\delta}, \enspace f_t(x)=f_t^{1,\delta}(x)+f_t^{2,\delta}(x)$, where
\begin{equation*}
\begin{aligned}
    b_{n,t,j}^{1,\delta}&:=\sum_{\substack{1\le i\le \lfloor nt \rfloor\\ |i-j|\le\delta n}}\alpha_{i-j}, \quad b_{n,t,j}^{2,\delta}:=\sum_{\substack{1\le i\le \lfloor nt \rfloor\\ |i-j|>\delta n}}\alpha_{i-j}, \\
    f_t^{1,\delta}(x)&:=\int_0^t\bigl[c_+(u-x)_+^{-\beta}+c_-(x-u)_+^{-\beta}\bigr]\mathbf1_{\{|u-x|\le\delta\}}du,\\
    f_t^{2,\delta}(x)&:=\int_0^t\bigl[c_+(u-x)_+^{-\beta}+c_-(x-u)_+^{-\beta}\bigr]\mathbf1_{\{|u-x|>\delta\}}du.
\end{aligned}
\end{equation*}
By \eqref{eq:alpha-global-bound} and $\beta<1$, for all sufficiently large $n$,
\begin{equation*}
\begin{aligned}
    |b_{n,t,j}^{1,\delta}|&\leq \sum_{\substack{1\le i\le \lfloor nt \rfloor\\ |i-j|\le\delta n}}|\alpha_{i-j}|\leq C\sum_{ |i-j|\le\delta n}(1+|i-j|)^{-\beta}\le  C(\delta n)^{1-\beta}, \\
    |f_t^{1,\delta}(x)|
 &\le
 |c_+|\int_x^{x+\delta}(u-x)^{-\beta}\,du
 +|c_-|\int_{x-\delta}^{x}(x-u)^{-\beta}\,du=\frac{|c_+|+|c_-|}{1-\beta}\,\delta^{1-\beta}.
\end{aligned}
\end{equation*}
Consequently,
\begin{equation}
\sup_j\bigl|n^{\beta-1}b_{n,t,j}^{1,\delta}\bigr|
\le C\delta^{1-\beta},
\qquad
\sup_x|f_t^{1,\delta}(x)|\le C\delta^{1-\beta}.
\label{eq:local-near-bounds}
\end{equation}

For the far part, we now establish
\begin{equation}
    \sup_{|j|\leq Rn}|n^{\beta-1}b_{n,t,j}^{2,\delta}-f_t^{2,\delta}(j/n)|\to 0.
\label{eq:far-part}
\end{equation}
Define
\begin{equation*}
\eta_r^+:=r^\beta\alpha_r-c_+,\qquad
\eta_r^-:=r^\beta\alpha_{-r}-c_-,\qquad r\ge1.
\end{equation*}
Then $\eta_r^\pm\to0$, and for fixed $\delta>0$,
\begin{equation*}
\varepsilon_{n,\delta}
:=\max\left\{
\sup_{r>\delta n}|\eta_r^+|,
\sup_{r>\delta n}|\eta_r^-|
\right\}\longrightarrow0.
\end{equation*}
On $i-j>\delta n$, we can write $n^\beta \alpha_{i-j}=(i/n-j/n)^{-\beta}(c_++\eta_{i-j}^+)$, so that
\begin{equation*}
    |n^\beta\alpha_{i-j}-c_+(i/n-j/n)^{-\beta}|\leq (i/n-j/n)^{-\beta} \varepsilon_{n,\delta} \le \delta^{-\beta}\varepsilon_{n,\delta}.
\end{equation*}
For $j-i>\delta n$, the negative-side analogue gives
\begin{equation*}
    |n^\beta\alpha_{i-j}-c_-(j/n-i/n)^{-\beta}|\leq \delta^{-\beta}\varepsilon_{n,\delta}.
\end{equation*}
Set
\begin{equation*}
K_\delta(u,x)
:=\bigl[c_+(u-x)_+^{-\beta}+c_-(x-u)_+^{-\beta}\bigr]
\mathbf1_{\{|u-x|>\delta\}}.
\end{equation*}
It follows that 
\begin{align*}
    &\left|n^{\beta -1}b_{n,t,j}^{2,\delta}-\frac{1}{n}\sum_{i=1}^{\lfloor nt \rfloor}K_\delta(i/n, j/n)\right|\nonumber\\
    =&\frac{1}{n}\bigg|\sum_{\substack{1\le i\le \lfloor nt \rfloor\\ |i-j|>\delta n}}\big[n^\beta \alpha_{i-j}-(c_+(i/n-j/n)_+^{-\beta}+c_-(j/n-i/n)_+^{-\beta})\big]\bigg|\nonumber\\
    \leq&\frac{1}{n}nt\delta^{-\beta}\varepsilon_{n,\delta}=t\delta^{-\beta}\varepsilon_{n,\delta}.
\end{align*}
Because the estimate is uniform in $j$, we have
\begin{equation}
    \sup_j \left|n^{\beta -1}b_{n,t,j}^{2,\delta}-\frac{1}{n}\sum_{i=1}^{\lfloor nt \rfloor}K_\delta(i/n, j/n)\right|\le t\delta^{-\beta}\varepsilon_{n,\delta}\to0.
\end{equation}
Therefore, to complete the proof of \eqref{eq:far-part}, it suffices to verify
\begin{equation*}
    \sup_{|j|\leq Rn}\left|\frac{1}{n}\sum_{i=1}^{\lfloor nt \rfloor}K_\delta(i/n, j/n)-\int_0^tK_\delta(u,j/n)du\right|\to0.
\end{equation*}
For fixed $x$, the map
$u\mapsto K_\delta(u,x)$ is bounded by
$M_\delta:=(|c_+|+|c_-|)\delta^{-\beta}$ and is
$L_\delta$-Lipschitz, on each interval separated by
$x-\delta$ and $x+\delta$, with
$L_\delta:=\beta(|c_+|+|c_-|)\delta^{-\beta-1}$. 
For the right term, we have the approximation 
\begin{equation}
 \left|\int_0^tK_\delta(u,j/n)du-\frac{1}{n}\int_0^{\lfloor nt\rfloor}K_\delta(u/n,j/n)du\right|=\left|\frac{1}{n}\int_{\lfloor nt\rfloor}^{nt}K_\delta(u/n,j/n)du\right|\leq \frac{M_\delta}{n}.
 \label{eq:far-approx1}
\end{equation}
Moreover, since $u\mapsto K_\delta(u,x)$ is $L_\delta$-Lipschitz on at least $\lfloor nt\rfloor-2$ intervals, it follows that
\begin{align}
    &\left|\frac{1}{n}\sum_{i=1}^{\lfloor nt \rfloor}K_\delta(i/n, j/n)-\frac{1}{n}\int_0^{\lfloor nt\rfloor}K_\delta(u/n,j/n)du\right|\nonumber\\
    \leq &\left|\frac{1}{n}\sum_{i=1}^{\lfloor nt \rfloor}\big[K_\delta(i/n,j/n)-\int_{i-1}^i K_\delta(u/n,j/n)du\big]\right|\nonumber\\
    \leq & \frac{1}{n}\sum_{i=1}^{\lfloor nt \rfloor} \int_{i-1}^i |K_\delta(i/n, j/n)-K_\delta(u/n,j/n)|du\leq  \frac{tL_\delta}{n}+\frac{4M_\delta}{n}.
\label{eq:far-approx2}
\end{align}
Combining \eqref{eq:far-approx1} and \eqref{eq:far-approx2}, we obtain
\begin{align*}
    &\sup_{|j|\leq Rn}\left|\frac{1}{n}\sum_{i=1}^{\lfloor nt \rfloor}K_\delta(i/n, j/n)-\int_0^tK_\delta(u,j/n)du\right|\leq \frac{tL_\delta}{n}+\frac{5M_\delta}{n}\to 0.
\end{align*}
Thus \eqref{eq:far-part} holds. Together with \eqref{eq:local-near-bounds} and the fact that $1/2<\beta<1$, this implies
\begin{align}
    &\limsup_{n\to\infty}\sup_{|j|\le Rn}\bigl|n^{\beta-1}b_{n,t,j}-f_t(j/n)\bigr|\nonumber\\
    \le & \limsup_{n\to\infty}\left(\sup_{|j|\leq Rn}|n^{\beta-1}b_{n,t,j}^{2,\delta}-f_t^{2,\delta}(j/n)| + \sup_{|j|\leq Rn}|n^{\beta-1}b_{n,t,j}^{1,\delta}| + \sup_{|j|\leq Rn}|f_t^{1,\delta}(j/n)|\right)\nonumber\\
    \le & 2C\delta^{1-\beta}.
\end{align}
Letting $\delta\searrow0$ yields \eqref{eq:coefapprox}.
\end{proof}

We next use the i.i.d.\ invariance principle~\cite{Zhang2023} to establish
convergence of weighted sums over a fixed bounded interval.
\begin{lemma}\label{lem:main-approx}
    Under Assumption~\ref{ass:innovations}, let $\bm{f}:\mathbb{R}\mapsto \mathbb{R}^d$ be a continuous function and let $R_1<R_2$ be fixed integers. We denote $\bm{f}(x)=(f_{t_1}(x),\cdots,f_{t_d}(x))^\top$ and $\bm{T}_n=\sum_{R_1n\leq j\leq R_2n}\bm{b}_{n,j}\varepsilon_j$, where
    \begin{equation}
        \bm{b}_{n,j}=n^{-H}(b_{n,t_1, j},\cdots, b_{n,t_d,j})^\top.
    \end{equation}
    Suppose 
    \begin{equation}
        \sup_{R_1n\leq j\leq R_2n}\left\Vert\sqrt{n}\bm{b}_{n,j}-\bm{f}(j/n)\right\Vert_2\to0\quad \text{as }n\to\infty.
        \label{eq:kernel approx}
    \end{equation}
    Then for any $\varphi\in C_{b,Lip}(\mathbb{R}^d)$,
    \begin{equation}
       \lim_{n\to\infty} \wE[\varphi(\bm{T}_n)]=\sup_{\theta\in\mathcal{A}_{[R_1,R_2]}}E\left[\varphi\left(\int_{R_1}^{R_2}\bm{f}(x)\theta_xdW_x\right)\right].
    \end{equation}
\end{lemma}

\begin{proof}
    By translating $j$ and $x$, we may assume without loss of generality that $R_1 = 0$ and $R=R_2>0$. 
    Denote $\bm{s}_{n,j}:=\sqrt{n}\bm{b}_{n,j}-\bm{f}(j/n)$. Applying Corollary \ref{cor:rosenthal}, it follows that 
    \begin{align}
        &\left|\wE\bigg[\varphi\bigg(\sum_{j=0}^{Rn}\bm{b}_{n,j}\varepsilon_j\bigg)\bigg]-
        \wE\bigg[\varphi\bigg(\frac{1}{\sqrt{n}}\sum_{j=1}^{Rn}\bm{f}(j/n)\varepsilon_j\bigg)\bigg]\right|
        \le  \frac{C}{\sqrt{n}}\bE\bigg[\big\Vert\sum_{j=1}^{Rn}\bm{s}_{n,j}\varepsilon_j\big\Vert_2^2\bigg]^{1/2}+C\bE\big[\|\bm b_{n,0}\varepsilon_0\|_2^2\big]^{1/2}\nonumber\\
        \le & \frac{C}{\sqrt{n}}\bigg\{\sum_{j=1}^{Rn}\wE[\Vert \bm{s}_{n,j}\Vert_2^2]\overline{\sigma}^2\bigg\}^{1/2} + \frac{C\usigma}{\sqrt n}\big(\|\bm{f}(0)\|_2 +o(1)\big)\nonumber\\
        \le &  C\usigma\bigg(\sqrt{R}\sup_{0\le j\le Rn}\Vert \bm{s}_{n,j}\Vert_2 + \frac{1}{\sqrt{n}}(\|\bm{f}(0)\|_2 +o(1))\bigg) \to0.\label{eq:approx3}
    \end{align}
   Thus, it suffices to show 
    \begin{equation}
    \label{eq:target1}
        \wE\bigg[\varphi\bigg(\frac{1}{\sqrt{n}}\sum_{j=1}^{Rn}\bm{f}(j/n)\varepsilon_j\bigg)\bigg] \longrightarrow \sup_{\theta\in\mathcal{A}_{[0,R]}}E\left[\varphi\left(\int_{0}^{R}\bm{f}(x)\theta_xdW_x\right)\right].
    \end{equation}

    Without loss of generality, we can assume $\bm{f}\in C^\infty[0,R]$. Otherwise, for any continuous function $\bm{f}$ and $\delta>0$, there exists $\hat{\bm{f}}\in C^\infty[0,R]$ such that $\sup_{x\in[0,R]}\Vert \bm{f}(x)-\hat{\bm{f}}(x)\Vert_2<\delta$. By the same argument as in \eqref{eq:approx3}, we have
    \begin{align}
        &\left|\wE\bigg[\varphi\bigg(\frac{1}{\sqrt{n}}\sum_{j=1}^{Rn}\bm{f}(j/n)\varepsilon_j\bigg)\bigg]- \wE\bigg[\varphi\bigg(\frac{1}{\sqrt{n}}\sum_{j=1}^{Rn}\hat{\bm{f}}(j/n)\varepsilon_j\bigg)\bigg]\right|\nonumber\\
        \le &\frac C{\sqrt{n}}\bE\bigg[\Vert\sum_{j=1}^{Rn}(\bm{f}(j/n)-\hat{\bm{f}}(j/n))\varepsilon_j\Vert_2^2\bigg]^{1/2}\nonumber\\
        \le &C\usigma\sqrt{R}\sup_{0\le x\le R}\Vert\bm{f}(x)-\hat{\bm{f}}(x) \Vert_2  
        \le C \delta. \label{eq:approx4}
    \end{align}
    For the right side of \eqref{eq:target1}, it follows from It\^o isometry that
    \begin{align}
        &\sup_{\theta_x\in\mathcal{A}_{[0,R]}}\bigg|E\bigg[\varphi\bigg(\int_{0}^R\bm{f}(x)\theta_xdW_x\bigg)\bigg]-E\bigg[\varphi\bigg(\int_{0}^R\hat{\bm{f}}(x)\theta_xdW_x\bigg)\bigg]\bigg| \nonumber\\
        \le &C \sup_{\theta_x\in\mathcal{A}_{[0,R]}}\left( E\int_0^R\|f(x)-\widehat f(x)\|_2^2\theta_x^2\,dx\right)^{1/2}\nonumber\\
        = & C \usigma \left(\int_0^R\|f(x)-\widehat f(x)\|_2^2\,dx\right)^{1/2}\le C\usigma\sqrt{R}\delta.\label{eq:approx5}
    \end{align}

    We decompose
    \begin{align*}
        \frac{1}{\sqrt{n}}\sum_{j=1}^{Rn}\bm{f}(j/n)\varepsilon_j &=\frac{1}{\sqrt{n}}\sum_{j=1}^{Rn}\varepsilon_j\int_{j-1}^j\bm{f}(x/n)dx+\frac{1}{\sqrt{n}}\sum_{j=1}^{Rn}\varepsilon_j\int_{j-1}^j\big(\bm{f}(j/n)-\bm{f}(x/n)\big)dx \\
        &=:I_{n}+J_{n}.
    \end{align*}
    For $J_{n}$, note that 
    \begin{align*}
        \Vert \bm{s}'_{n,j} \Vert_2^2:=\bigg\Vert \int_{j-1}^j \big(\bm{f}(j/n)-\bm{f}(x/n)\big)dx\bigg\Vert_2^2\leq\int_{j-1}^j \big\Vert\bm{f}(j/n)-\bm{f}(x/n)\big\Vert_2^2dx\le \frac{\Vert \nabla \bm{f}\Vert_\infty^2}{n^2}.
    \end{align*}
    Applying Corollary \ref{cor:rosenthal}   again, it yields
    \begin{align}
        \bE[\Vert J_{n}\Vert_2^2]&\le\frac{C}{n} \sum_{j=1}^{Rn}\Vert \bm{s}'_{n,j}\Vert_2^2\overline{\sigma}^2  
       \le \frac{C}{n^2} R\usigma^2 \Vert\nabla \bm{f}\Vert_\infty^2 \to 0.\label{eq:approx6}
    \end{align}
    
    Combining \eqref{eq:approx4}--\eqref{eq:approx6}, our target \eqref{eq:target1} is reduced to 
    \begin{equation}
    \label{eq:target2}
        \wE\big[\varphi(I_n)\big] \longrightarrow \sup_{\theta\in\mathcal{A}_{[0,R]}}E\left[\varphi\left(\int_{0}^{R}\bm{f}(x)\theta_xdW_x\right)\right],\quad \bm{f}\in C^\infty[0,R].
    \end{equation}
    
    Denote $Q_0=0, Q_n=\sum_{j=1}^n\varepsilon_j$ and $Q_x=Q_n+(x-n)(Q_{n+1}-Q_n)$ for $x\in[n,n+1]$. Noting that for $x\in[n-1,n]$, $dQ_x=(Q_{n}-Q_{n-1})dx=\varepsilon_{n}dx$, it follows from direct integration by parts that
    \begin{align}
        I_n
         &=\frac1{\sqrt n}\int_0^{Rn}\bm f(x/n)\,dQ_x\nonumber\\
         &=\frac{Q_{Rn}}{\sqrt n}f(R)
           -\frac1{n\sqrt n}\int_0^{Rn}Q_x \nabla\bm f(x/n)\,dx\nonumber\\
         &=\frac{Q_{Rn}}{\sqrt n}f(R)
           -\int_0^R\frac{Q_{nu}}{\sqrt n}\nabla\bm f(u)\,du.
    \end{align}
    By the functional central limit theorem of i.i.d. random variables under sublinear expectations (see Theorem 3.6 in \cite{Zhang2023}), we have $\frac{Q_{nx}}{\sqrt{n}}\Rightarrow B_x^G$, where $B_x^G$ is a $G$-Brownian Motion on $\widetilde{\mathbb{E}}$ with $\widetilde{\mathbb{E}}[(B_1^G)^2]=\bE[\varepsilon_1^2]=\usigma^2$ and $\widetilde{\mathcal{E}}[(B_1^G)^2]=\be[\varepsilon_1^2]=\lsigma^2$. 

    Note that the mapping $\phi:W \longmapsto W_R\bm{f}(R)-\int_0^RW_x\nabla\bm{f}(x)dx$ is a continuous linear function of $W$. Hence $\varphi\circ\phi\in C_{b,Lip}(C[0,R])$. It yields that
    \begin{equation}\label{eq:approx8}
        \wE[\varphi(I_n)]=\wE\bigg[\varphi\circ\phi\big(\frac{Q_{n\cdot}}{\sqrt{n}}\big)\bigg]\longrightarrow \widetilde{\E}[\varphi\circ\phi(B^G)]= \widetilde{\mathbb{E}}\bigg[\varphi\bigg(B_R^G\bm{f}(R)-\int_0^RB_x^G\nabla\bm{f}(x)dx\bigg)\bigg].
    \end{equation}
    It yields from the representation of $G$-Brownian Motion \cite{DenisHuPeng2011} that 
    for any $\varphi\in C_{b,Lip}(\mathbb R^d)$ and $t_1,t_2,\ldots, t_d\in [0,R]$,
   \begin{align}\label{eq:distidenticalBG} \widetilde{\mathbb{E}}[\varphi(B_{t_1}^G, \ldots, B_{t_d}^G)]=&\sup_{P\in \mathcal{Q}}E_P[\varphi(B_{t_1}^G, \ldots, B_{t_d}^G)]\nonumber\\
   =&\sup_{\theta\in\mathcal{A}_{[0,R]}}E\bigg[\varphi\bigg(\int_0^{t_1} \theta_x d W_x, \ldots,\int_0^{t_d}\theta_x dW_x\bigg)\bigg], 
   \end{align}
     where $W_x$ is a Brownian Motion in a classical probability space. Let $P_{\theta}$ be the law of the process $W_t^{\theta}=\int_0^t \theta_s\;d W_s$. Notice that
     $$ \sup_{Q\in \mathcal{Q}}E_Q[|B_t^G-B_s^G|^4]=\widetilde{\mathbb E}[|B_t^G-B_s^G|^4]\le 3 \usigma^4(t-s)^2, $$
   and
   $$ \sup_{\theta\in \mathcal{A}}E[|W_t^{\theta}-W_s^{\theta}|^4]=\widetilde{\mathbb E}[|B_t^G-B_s^G|^4]\le 3 \usigma^4(t-s)^2. $$
  It follows that   both $\{QB^{-1}, Q\in \mathcal{Q}\}$ and $\{P_{\theta}, \theta\in \mathcal{A}[0,R]\}$ are tight on $C([0,R])$ (cf. Proposition 6.2.15 of Peng \cite{Peng2019} and its proof), and \eqref{eq:distidenticalBG} implies that 
    for any bounded Lipschitz function $h:C[0,R]\longmapsto \mathbb{R}$, 
    \begin{equation}\label{eq:representation}
        \widetilde{\mathbb E}
         \!\left[h\!\left((B_t^G)_{t\in[0,R]}\right)\right]
        =
        \sup_{\theta\in\mathcal A_{[0,R]}}
        E\!\left[
        h\!\left(
        \left(\int_0^t\theta_x\,dW_x\right)_{t\in[0,R]}
        \right)\right].
    \end{equation}
   cf. Lemma 3.2 of Zhang \cite{Zhang2026}. See Proposition \ref{prop:path-tightness-1} below. We can verify that
    \begin{align}
        \phi\bigg(\int_0^t\theta_xdWx\bigg)&=\bm{f}(R)\int_0^R\theta_xdW_x-\int_0^R\int_{0}^x\theta_sdW_s\, \nabla\bm{f}(x)dx\nonumber\\
        &=\bm{f}(R)\int_0^R\theta_xdW_x-\int_0^R\int_s^R\nabla\bm{f}(x)dx\, \theta_sdW_s\nonumber\\
        &=\bm{f}(R)\int_0^R\theta_xdW_x-\int_0^R\big(\bm{f}(R)-\bm{f}(s)\big)\theta_sdW_s\nonumber\\
        &=\int_0^R \bm{f}(x)\theta_xdW_x. \label{eq:integral_simplify}
    \end{align}
    Substituting \eqref{eq:representation} and \eqref{eq:integral_simplify} into the right side of \eqref{eq:approx8} obtains
    \begin{align*}
        &\widetilde{\mathbb{E}}\bigg[\varphi\bigg(B_R^G\bm{f}(R)-\int_0^RB_x^G\nabla\bm{f}(x)dx\bigg)\bigg]=\widetilde{\mathbb{E}}\big[\varphi \circ\phi(B_t^G)\big]\\
        =&\sup_{\theta\in\mathcal{A}_{[0,R]}} E\bigg[\varphi\circ\phi\bigg(\int_0^t\theta_xdWx\bigg)\bigg]=\sup_{\theta\in\mathcal{A}_{[0,R]}} E\bigg[\varphi\bigg(\int_0^R\bm{f}(x)\theta_xdWx\bigg)\bigg].
    \end{align*}
\end{proof}
Combining the preceding lemmas with a tail estimate, we obtain
finite-dimensional convergence of the linearly interpolated process.
\begin{theorem}\label{thm:finite-dimensional-convergence}
Under Assumptions~\ref{ass:innovations} and~\ref{ass:coeff},
let $\widetilde X_n$ be the linear interpolation
defined in \eqref{eq:defX}.
For every $d\ge1$, $0\le t_1<\cdots<t_d\le1$, and
$\varphi\in C_{b,\mathrm{Lip}}(\mathbb R^d)$,
$$
 \lim_{n\to\infty}
 \widehat{\mathbb E}
 \bigl[\varphi(\widetilde X_n(t_1),\ldots,
              \widetilde X_n(t_d))\bigr]=\sup_{\theta\in\mathcal A}
 E\left[
 \varphi\left(\int_{\mathbb R}\bm{f}(x)\theta_x\,dW_x\right)
 \right],
$$
where $\bm{f}(x):=(f_{t_1}(x),\ldots,f_{t_d}(x))^\top$
and $\mathcal A$ is the whole-line predictable control
class in \eqref{eq:whole-line-controls}.
\end{theorem}

\begin{proof}
Set $\bm{T}_n:=\sum_{j\in\mathbb Z}\bm{b}_{n,j}\varepsilon_j$. By \eqref{eq:wp-partial-sums}, $ \bm{T}_n=n^{-H}(S_{\lfloor nt_1\rfloor},\ldots,S_{\lfloor nt_d\rfloor})^\top.$ 
Let
\[
 L(\varphi)
 :=
 \sup_{\theta\in\mathcal A}
 \mathbb E\left[
 \varphi\left(\int_{\mathbb R}f(x)\theta_x\,dW_x\right)
 \right].
\]
We first prove that
$\widehat{\mathbb E}[\varphi(\bm{T}_n)]\to L(\varphi)$.

For an integer $R\ge2$, set $\bm{T}_{n,R}:=\sum_{|j|\le Rn}\bm{b}_{n,j}\varepsilon_j.$
Since $1/2-H=\beta-1$, Lemma~\ref{lem:local-profile}
gives
\begin{align*}
    \sup_{|j|\le Rn}
 \|\sqrt n\,\bm{b}_{n,j}-\bm{f}(j/n)\|_2\le\bigg(\sum_{r=1}^d\big(\sup_{|j|\le Rn}|n^{\beta-1}b_{n,t_r,j}-f_{t_r}(j/n)|\big)^2\bigg)^{1/2}\longrightarrow0.
\end{align*}
Lemma~\ref{lem:main-approx} therefore yields
\begin{equation}\label{eq:fixed-R-fdd}
 \widehat{\mathbb E}[\varphi(\bm{T}_{n,R})]
 \longrightarrow
 L_R(\varphi)
 :=
 \sup_{\theta\in\mathcal A_{[-R,R]}}
 E\left[
 \varphi\left(\int_{-R}^{R}\bm{f}(x)\theta_x\,dW_x\right)
 \right].
\end{equation}

We next estimate the contribution from $|j|>Rn$.
The coefficient assumption implies
$|\alpha_k|\le C(1+|k|)^{-\beta}$.
For $|j|>Rn$ and $1\le i\le n$, we have
$|i-j|\ge |j|/2$, and hence
\[
 |b_{n,t,j}|
 \le
 \sum_{i=1}^{\lfloor nt\rfloor}|\alpha_{i-j}|
 \le Cn|j|^{-\beta},
 \qquad 0\le t\le1.
\]
Consequently,
\begin{equation}\label{eq:coefficient-tail}
 \sum_{|j|>Rn}\|\bm{b}_{n,j}\|_2^2
 \le Cn^{2-2H}\sum_{|j|>Rn}|j|^{-2\beta}
 \le CR^{1-2\beta},
\end{equation}
where the last inequality uses $2H=3-2\beta$.
The constants only depend on
$d,t_1,\ldots,t_d$, but not on $n$ or $R$.

By Theorem~\ref{thm:wp-weighted-series} and the
Lipschitz property of $\varphi$,

\begin{equation}\label{eq:discrete-tail}
 \left|
 \widehat{\mathbb E}[\varphi(\bm{T}_n)]
 -\widehat{\mathbb E}[\varphi(\bm{T}_{n,R})]
 \right|\le
 C\left(
 \widehat{\mathbb E}[\|\bm{T}_n-\bm{T}_{n,R}\|_2^2]
 \right)^{1/2}\le
 C\bar\sigma
 \left(\sum_{|j|>Rn}\|\bm{b}_{n,j}\|_2^2\right)^{1/2}
 \le C R^{1/2-\beta}.
\end{equation}

Thus the discrete tail is negligible uniformly in $n$.

We next remove the truncation in the limiting integral.  
For $|x| > R \ge 2$ and $u \in [0, 1]$, $|u - x| \ge |x|/2$, 
from the definition of $f_t$,
$\|\bm{f}(x)\|_2\le C|x|^{-\beta}$ for $|x|\ge2$, so
\begin{equation}\label{eq:f-tail}
 \int_{|x|>R}\|\bm{f}(x)\|_2^2\,dx
 \le CR^{1-2\beta}.
\end{equation}
To compare $L_R(\varphi)$ with $L(\varphi)$,
we first show that
\begin{equation}\label{eq:control-class-identification}
 L_R(\varphi)
 =
 \sup_{\theta\in\mathcal A}
  E\left[
 \varphi\left(\int_{-R}^{R}\bm{f}(x)\theta_x\,dW_x\right)
 \right].
\end{equation}
Identifying each finite-interval control with its
extension taking the constant value $\underline\sigma$
outside $[-R,R]$, we have
$\mathcal A_{[-R,R]}\subset\mathcal A$.
Hence
\[
 L_R(\varphi)
 \le
 \sup_{\theta\in\mathcal A}
  E\left[
 \varphi\left(
 \int_{-R}^{R}\bm{f}(x)\theta_x\,dW_x
 \right)\right].
\]

Conversely, let
\[
 \mathcal G_R
 :=\sigma\{W_v-W_u:u<v\le -R\}.
\]
For any $\theta\in\mathcal A$, conditioning on
$\mathcal G_R$ fixes the past information used by
the control, while the increments on $[-R,R]$
retain their Brownian law.
The resulting strategy is therefore an admissible
finite-interval control, so
\[
  E\left[
 \varphi\left(
 \int_{-R}^{R}\bm{f}(x)\theta_x\,dW_x
 \right)
 \,\middle|\,\mathcal G_R
 \right]
 \le L_R(\varphi)
 \qquad\text{a.s.}
\]
Taking expectations and then the supremum over
$\theta\in\mathcal A$ yields
\[
 \sup_{\theta\in\mathcal A}
  E\left[
 \varphi\left(
 \int_{-R}^{R}\bm{f}(x)\theta_x\,dW_x
 \right)\right]
 \le L_R(\varphi).
\]
Together, the two inequalities establish
\eqref{eq:control-class-identification}.

It follows from \eqref{eq:control-class-identification},
It\^o's isometry, $|\theta_x|\le\bar\sigma$, and \eqref{eq:f-tail} that

\begin{align}
 |L_R(\varphi)-L(\varphi)|
 &\le
 C\sup_{\theta\in\mathcal A}
 \left( E\left[
 \left\|
 \int_{|x|>R}\bm{f}(x)\theta_x\,dW_x
 \right\|_2^2
 \right]
 \right)^{1/2}\nonumber\\
 &\le
 C\bar\sigma 
 \left(
 \int_{|x|>R}\|\bm{f}(x)\|_2^2\,dx
 \right)^{1/2}\le C R^{1/2-\beta}.\label{eq:integral-tail}
\end{align}

Combining the two tail estimates \eqref{eq:discrete-tail} and \eqref{eq:integral-tail} with
\eqref{eq:fixed-R-fdd} and letting $R\to\infty$, we obtain from $\beta>1/2$ that
\[
 \limsup_{n\to\infty}
 \left|
 \widehat{\mathbb E}[\varphi(\bm{T}_n)]-L(\varphi)
 \right|
 \le C_\varphi R^{1/2-\beta}\longrightarrow 0.
\]

It remains to account for the linear interpolation.
Put
\[
 \widetilde{\bm X}_n
 :=
 (\widetilde X_n(t_1),\ldots,
  \widetilde X_n(t_d))^\top,
 \qquad
 \lambda_{n,r}:=nt_r-\lfloor nt_r\rfloor.
\]
The $r$-th component of $\widetilde{\bm X}_n-\bm{T}_n$ is
\[
 n^{-H}\lambda_{n,r}Y_{\lfloor nt_r\rfloor+1}.
\]
Since $0\le\lambda_{n,r}\le1$,
Corollary~\ref{cor:wp-linear-process} gives
\[
\begin{aligned}
 \widehat{\mathbb E}
 [\|\widetilde{\bm X}_n-\bm{T}_n\|_2^2]
 \le
 n^{-2H}\sum_{r=1}^d
 \widehat{\mathbb E}
 [|Y_{\lfloor nt_r\rfloor+1}|^2]\le
 d\bar\sigma^2 n^{-2H}
 \sum_{j\in\mathbb Z}\alpha_j^2
 \longrightarrow0.
\end{aligned}
\]
Therefore,
\[
 \left|
 \widehat{\mathbb E}[\varphi(\widetilde{\bm X}_n)]
 -\widehat{\mathbb E}[\varphi(\bm{T}_n)]
 \right|
 \le
 C\left(
 \widehat{\mathbb E}
 [\|\widetilde{\bm X}_n-\bm{T}_n\|_2^2]
 \right)^{1/2}
 \longrightarrow0.
\]
This completes the proof.
\end{proof}

\section{Tightness and the proof of the main result}
\label{sec:tightness}

We now establish tightness and complete the proof of
Theorem~\ref{thm:main}. The argument uses only the
second-moment estimate of
Theorem~\ref{thm:wp-weighted-series}.

\begin{lemma}\label{lem:increment-second-moment}
Under Assumptions~\ref{ass:innovations}
and~\ref{ass:coeff}, there exists a constant $C<\infty$
such that
\begin{equation}\label{eq:tight-increment}
  \widehat{\mathbb E}
  \bigl[
    |\widetilde X_n(t)-\widetilde X_n(s)|^2
  \bigr]
  \le C|t-s|^{2H},
  \qquad n\ge1,\quad s,t\in[0,1].
\end{equation}
\end{lemma}

\begin{proof}
Assumption~\ref{ass:coeff} implies
\[
  |\alpha_r|\le C(1+|r|)^{-\beta},
  \qquad r\in\mathbb Z.
\]
For an integer $m\ge1$, put
\[
  a_{m,j}:=\sum_{i=1}^{m}\alpha_{i-j}.
\]
We first show that
\begin{equation}\label{eq:tight-coefficient-bound}
  \sum_{j\in\mathbb Z}|a_{m,j}|^2
  \le C m^{3-2\beta}
  =C m^{2H}.
\end{equation}
If $|j|\le2m$, since $|i-j|\le 3m$, then
\[
  |a_{m,j}|
  \le C\sum_{|r|\le3m}(1+|r|)^{-\beta}
  \le C m^{1-\beta}.
\]
If $|j|>2m$, then $|i-j|\ge |j|/2$ for
$1\le i\le m$, and hence
\[
  |a_{m,j}|\le \sum_{i=1}^m |\alpha_{i-j}|\le Cm(1+|j|/2)^{-\beta}\le C m|j|^{-\beta}.
\]
Consequently,
\[
\begin{aligned}
  \sum_{j\in\mathbb Z}|a_{m,j}|^2\le
  C m^{2-2\beta}\#\{j:|j|\le2m\}
  +C m^2\sum_{|j|>2m}|j|^{-2\beta}\le C m^{3-2\beta}=Cm^{2H}.
\end{aligned}
\]

For integers $0\le \ell<k$, the series representation
of the partial sums gives
\[
  S_k-S_\ell
  =
  \sum_{j\in\mathbb Z}
  \left(\sum_{i=\ell+1}^{k}\alpha_{i-j}\right)
  \varepsilon_j
  \quad\text{in }L_{\mathcal P}^2.
\]
By \eqref{eq:limit-l2-norm}, and
\eqref{eq:tight-coefficient-bound},
\begin{equation}\label{eq:tight-grid-increment}
\begin{aligned}
  \widehat{\mathbb E}[|S_k-S_\ell|^2]\le
  \bar\sigma^2
  \sum_{j\in\mathbb Z}
  \left|\sum_{i=\ell+1}^{k}\alpha_{i-j}\right|^2\le C(k-\ell)^{2H}.
\end{aligned}
\end{equation}
In particular, $\sup_{i\ge1}\|Y_i\|_{2,\mathcal P}<\infty$.

It remains to pass from grid points to arbitrary times.
Fix $0\le s<t\le1$.
On each interval $[(k-1)/n,k/n]$, the interpolation
has slope $n^{1-H}Y_{k}$.

If $(k-1)/n\le s<t\le k/n$ for some integer $k$, since $(n|t-s|)^{1-H}\le 1$, then
\begin{align}\label{eq:interval-increment}
    \|\widetilde X_n(t)-\widetilde X_n(s)\|_{2,\mathcal P}\le n^{1-H}|t-s|\|Y_k\|_{2,\mathcal{P}}\le C (n|t-s|)^{1-H}|t-s|^H\le C|t-s|^H.
\end{align}

Otherwise, set $s_n:=\frac{\lceil ns\rceil}{n}$ and $t_n:=\frac{\lfloor nt\rfloor}{n}.$ Then $s\le s_n\le t_n\le t$. It follows from \eqref{eq:tight-grid-increment} and \eqref{eq:interval-increment} that
\begin{align*}
    \|\widetilde X_n(t)-\widetilde X_n(s)\|_{2,\mathcal P}&\le \|\widetilde X_n(t)-\widetilde X_n(t_n)\|_{2,\mathcal P}+ \|\widetilde X_n(t_n)-\widetilde X_n(s_n)\|_{2,\mathcal P}+\|\widetilde X_n(s_n)-\widetilde X_n(s)\|_{2,\mathcal P}\\
    &\le C|t-t_n|^H+C|t_n-s_n|^H+C|s_n-s|^H\le C|t-s|^H.
\end{align*}
Squaring proves \eqref{eq:tight-increment}.

\end{proof}

For $m\ge0$, let $\pi_m:C[0,1]\to C[0,1]$ be
the linear interpolation operator on the dyadic grid
$\{k2^{-m}:0\le k\le2^m\}$. For $x\in C[0,1]$, 
\begin{equation}
    (\pi_mx)(t):=(1-\lambda_t)x(k2^{-m})+\lambda_tx((k+1)2^{-m}),
\end{equation}
where $\lambda_t=2^mt-k\in[0,1]$ for $t\in[k2^{-m},(k+1)2^{-m}].$ 
\begin{proposition}\label{prop:path-tightness-1} Let $\mathcal{Q}$ be a family of probability measures on $C[0,1]$ and $X_t(x):=x(t)$ be the coordinate process on $C[0,1]$. Suppose $Q(X_0=0)=1$, $Q\in \mathcal{Q}$, and there exist positive constants $p\ge 1$, $\beta$, $\alpha$ such that 
\begin{equation}\label{eq:uniform-law-increments-1} \sup_{Q\in \mathcal{Q}}E_Q\big[|X_t-X_s|^p]\le \beta |t-s|^{\alpha+1}, \;\; s,t\in [0,1]. 
\end{equation}
Then the family $\mathcal{Q}$ is  tight on $C[0,1]$. Moreover, there exists $C<\infty$
such that, for every $m\ge0$,
\begin{equation}\label{eq:dyadic-path-approximation-1}
  \sup_{Q\in \mathcal{Q}}
  E_Q\bigl[
    \|X-\pi_m X\|_\infty
  \bigr]
  \le C2^{-m\alpha/p},
\end{equation}
and thus, for any $\Phi\in C_{b,Lip}(C[0,1])$, 
$$  \sup_{Q\in \mathcal{Q}}
  E_Q\bigl[\big|\Phi(X)-\Phi(\pi_m X)\big|
  \bigr]
  \le C \operatorname{Lip}(\Phi) 2^{-m\alpha/p}.$$
\end{proposition}

\begin{proof} 
For $x\in C[0,1]$ and $r\ge0$, define
\[
  D_r(x):=
  \max_{0\le k<2^r}
  \left|
    x((k+1)2^{-r})-x(k2^{-r})
  \right|.
\]
Then, \eqref{eq:uniform-law-increments-1} yields
\begin{align}\label{eq:dyadic-increment-bound-1}
  E_Q[D_r(X)^p]&=E_Q\big[\max_{0\le k< 2^r}|X_{(k+1)2^{-r}}-X_{k2^{-r}}|^p\big] \nonumber\\
  &\le
  \sum_{k=0}^{2^r-1}
  E_Q\bigl[
    |X_{(k+1)2^{-r}}-X_{k2^{-r}}|^p
  \bigr]\nonumber\\
  &\le \beta 2^r\cdot 2^{-r(\alpha+1)}= \beta 2^{-r\alpha}.
\end{align}

Choose $0<\gamma<\alpha/p$. For $L>0$, let
\[
  K_L:=
  \left\{
    x\in C[0,1]:
    x(0)=0,\quad
    D_r(x)\le L2^{-r\gamma}
    \text{ for every }r\ge0
  \right\}.
\]
We next show that  $K_L$ is compact and $Q(K_L^c)\le C_{\gamma}/L^p$.

Fix $x\in K_L$ and $0\le s<t\le1$.
Choose an integer $r\ge0$ such that $2^{-r-1}<t-s\le2^{-r}.$
For $u\in[0,1]$ and $q\ge r$, let $ u_q:=2^{-q}\lfloor2^q u\rfloor$ 
be the dyadic grid point immediately to the left of $u$.
Then $u_q\to u$ as $q\to\infty$.

Since $t-s\le2^{-r}$, the points $s_r$ and $t_r$
are equal or adjacent on the level-$r$ grid. So $t_r-s_r=0$ or $2^{-r}$.
Consequently,
\begin{equation}\label{eq:tight-coarse-increment-1}
  |x(t_r)-x(s_r)|
  \le D_r(x)
  \le L2^{-r\gamma}.
\end{equation}

Thus
\[
  |x(u_{q+1})-x(u_q)|
  \le D_{q+1}(x)
  \le L2^{-(q+1)\gamma}.
\]
For every $N>r$, summing these increments gives
\[
  |x(u_N)-x(u_r)|
  \le
  \sum_{q=r}^{N-1}|x(u_{q+1})-x(u_q)|
  \le
  L\sum_{q=r+1}^{N}2^{-q\gamma}.
\]
By continuity of $x$, letting $N\to\infty$ yields
\begin{equation}\label{eq:tight-grid-error-1}
  |x(u)-x(u_r)|
  \le
  L\sum_{q=r+1}^{\infty}2^{-q\gamma}
  =
  \frac{L2^{-r\gamma}}{2^\gamma-1},
  \qquad u=s,t.
\end{equation}

Combining \eqref{eq:tight-coarse-increment-1}
and \eqref{eq:tight-grid-error-1}, we obtain
\[
\begin{aligned}
  |x(t)-x(s)|
  &\le
  |x(t)-x(t_r)|
  +|x(t_r)-x(s_r)|
  +|x(s_r)-x(s)|\\
  &\le
  L\left(1+\frac{2}{2^\gamma-1}\right)
  2^{-r\gamma} 
 \le C_\gamma L|t-s|^\gamma,
\end{aligned}
\]
where the last inequality follows from
$2^{-r}<2(t-s)$.
The constant $C_\gamma$ is independent of
$x$, $s$, and $t$.

Since $x(0)=0$ for every $x\in K_L$, this estimate gives
\[
  \sup_{x\in K_L}\|x\|_\infty\le C_\gamma L.
\]
Define $\omega_{\delta}(x)=\sup_{|t-s|<\delta}|x(t)-x(s)|$. It also gives
\[
  \sup_{x\in K_L}\omega_\delta(x)
  \le C_\gamma L\delta^\gamma
  \longrightarrow0
  \qquad\text{as }\delta\downarrow0.
\]
Hence $K_L$ is uniformly bounded and equicontinuous.

Finally, $K_L$ is closed in the supremum norm.
Indeed, if $x_n\in K_L$ and $x_n\to x$ uniformly,
then $x(0)=0$ and, for each fixed $q$,
\[
  D_q(x)=\lim_{n\to\infty}D_q(x_n)
  \le L2^{-q\gamma}.
\]
Thus $x\in K_L$.
The Arzel\`a--Ascoli theorem now shows that $K_L$
is compact.

Since $Q(X_0=0)=1$, the union bound, Markov's inequality, and
\eqref{eq:dyadic-increment-bound-1} give
\[
\begin{aligned}
  Q(K_L^c)&\le Q\big(\bigcup_{r=0}^\infty\{D_r(X)>L2^{-r\gamma}\}\big)\le
  \sum_{r=0}^{\infty}
  Q\bigl(D_r(X)>L2^{-r\gamma}\bigr)\\
  &\le
  \frac{\beta}{L^p}
  \sum_{r=0}^{\infty}
  2^{-r(\alpha-p\gamma)}
  \le \frac{C_\gamma}{L^p}.
\end{aligned}
\]
This bound is uniform in $Q$ and proves tightness.

We finally establish
\eqref{eq:dyadic-path-approximation-1}.
For every continuous path $x$,
$\pi_r x\to x$ uniformly and
\[
  \|\pi_{r+1}x-\pi_r x\|_\infty
  \le D_{r+1}(x).
\]

Hence
\[
  \|x-\pi_m x\|_\infty
  \le \sum_{r=m+1}^{\infty}D_r(x).
\]
Using monotone convergence and
\eqref{eq:dyadic-increment-bound-1}, we obtain
\[
\begin{aligned}
  E_Q[\|X-\pi_m X\|_\infty]
  &\le
  \sum_{r=m+1}^{\infty}E_Q[D_r(X)]\le\sum_{r=m+1}^{\infty}E_Q[D_r(X)^p]^{1/p}\\
  &\le
  \beta^{1/p}\sum_{r=m+1}^{\infty}2^{-r\alpha/p}
  \le C2^{-m\alpha/p}.
\end{aligned}
\]
Taking the supremum over $Q\in\mathcal{Q}$ proves
\eqref{eq:dyadic-path-approximation-1}.
\end{proof}

\begin{proposition}\label{prop:path-tightness}
The families of probability measures
\begin{equation}\label{eq:laws}
      \bigl\{
    P\circ\widetilde X_n^{-1}:
    n\ge1,\ P\in\mathcal P
  \bigr\},
  \qquad
  \bigl\{
    P_W\circ(B^\theta)^{-1}:
    \theta\in\mathcal A
  \bigr\}
\end{equation}
are tight on $C[0,1]$. Moreover, there exists $C<\infty$
such that, for every $m\ge0$,
\begin{equation}\label{eq:dyadic-path-approximation}
  \sup_{n\ge1}
  \widehat{\mathbb E}
  \bigl[
    \|\widetilde X_n-\pi_m\widetilde X_n\|_\infty
  \bigr]+
  \sup_{\theta\in\mathcal A}
  E\bigl[
    \|B^\theta-\pi_m B^\theta\|_\infty
  \bigr]
  \le C2^{-m(H-1/2)}.
\end{equation}
\end{proposition}

\begin{proof}
Let $Q$ be any probability measure in either family in \eqref{eq:laws}, and let
$X_t(x):=x(t)$ be the coordinate process on $C[0,1]$.
Lemma~\ref{lem:increment-second-moment} gives the required
increment bound for the first family. For the second,
It\^o's isometry gives
\[
  E[|B^\theta_t-B^\theta_s|^2]
  \le
  \bar\sigma^2\|f_t-f_s\|_{L^2(\mathbb R)}^2
  =
  \bar\sigma^2\kappa^2|t-s|^{2H}.
\]
Therefore, for each $Q$, we have
\begin{equation}\label{eq:uniform-law-increments}
  Q(X_0=0)=1,
  \qquad
  E_Q[|X_t-X_s|^2]\le C|t-s|^{2H}.
\end{equation}
The conclusions follow from Proposition \ref{prop:path-tightness-1}. 
\end{proof}

\begin{proof}[Proof of Theorem~\ref{thm:main}]
We first assume   $\Phi\in C_{b,\mathrm{Lip}}(C[0,1]).$
For each fixed $m\ge0$, set
\[
  N:=2^m,\qquad t_k:=k/N,\quad k=0,\ldots,N.
\]
Let $g:\mathbb R^{N+1}\to C[0,1]$ map
$z=(z_0,\ldots,z_N)$ to the linear interpolation path taking
the value $z_k$ at $t_k$. Then
\[
  \pi_m x=g(x(t_0),\ldots,x(t_N)),
  \qquad x\in C[0,1].
\]
For any $z,w\in\mathbb R^{N+1}$, linear interpolation
gives
\[
  \|g(z)-g(w)\|_\infty
  \le\max_{0\le k\le N}|z_k-w_k|
  \le\|z-w\|_2.
\]
Consequently,
\[
  |(\Phi\circ g)(z)-(\Phi\circ g)(w)|
  \le C\|z-w\|_2,
  \qquad
  |(\Phi\circ g)(z)|\le\|\Phi\|_\infty.
\]
Thus $\Phi\circ g\in C_{b,\mathrm{Lip}}(\mathbb R^{N+1})$.
Theorem~\ref{thm:finite-dimensional-convergence}
therefore gives
\begin{align*}
      \lim_{n\to\infty}
  \widehat{\mathbb E}
  [\Phi(\pi_m\widetilde X_n)]
  &=  \lim_{n\to\infty}
  \widehat{\mathbb E}
  [(\Phi\circ g)(\widetilde{X}_{n}(t_0),\cdots,\widetilde{X}_{n}(t_N))]\\
  &= \sup_{\theta\in\mathcal A}
  E[(\Phi\circ g ) (B^\theta_{t_0},\cdots, B^\theta_{t_N})] \\
  &=\sup_{\theta\in\mathcal A}E[\Phi(\pi_m B^\theta)]
  =
  \widetilde{\mathbb E}[\Phi\circ\pi_m].
\end{align*}

On the other hand,
\eqref{eq:dyadic-path-approximation} implies
\[
\begin{aligned}
  \left|
    \widehat{\mathbb E}[\Phi(\widetilde X_n)]
    -
    \widehat{\mathbb E}[\Phi(\pi_m\widetilde X_n)]
  \right|\le
 C\,
  \widehat{\mathbb E}
  [\|\widetilde X_n-\pi_m\widetilde X_n\|_\infty]
  \le
  C\,2^{-m(H-1/2)}
\end{aligned}
\]
uniformly in $n$. Similarly,
\[
  \left|
    \widetilde{\mathbb E}[\Phi]
    -
    \widetilde{\mathbb E}[\Phi\circ\pi_m]
  \right|
  \le
  C\,2^{-m(H-1/2)}.
\]
Combining these estimates and first letting $n\to\infty$,
we obtain
\[
  \limsup_{n\to\infty}
  \left|
    \widehat{\mathbb E}[\Phi(\widetilde X_n)]
    -
    \widetilde{\mathbb E}[\Phi]
  \right|
  \le
  C\,2^{-m(H-1/2)}.
\]
Since $H>\frac12$, the right-hand side tends to zero
as $m\to\infty$. 

Now, fix $\Phi\in C_b(C[0,1])$. Without loss of generality we assume $\|\Phi\|=\sup_{\bm x\in C[0,1]}|\Phi(\bm x)|<1$. Then there exists a sequence of functions $\Phi_k\in C_{b,Lip}(C[0,1])$ with $\|\Phi_k\|\le 1$ such that
$$ \|\Phi-\Phi_k\|_K=:\sup_{\bm x\in K}|\Phi(\bm x)-\Phi_k(\bm x)|\to 0 \text{ for any compact set } K \subset C[0,1]. $$
Then
\begin{align*}
&\left|
    \widehat{\mathbb E}[\Phi(\widetilde X_n)]
    -
    \widetilde{\mathbb E}[\Phi]
  \right|\le \left|
    \widehat{\mathbb E}[\Phi_k(\widetilde X_n)]
    -
    \widetilde{\mathbb E}[\Phi_k]
  \right|+\left|
    \widehat{\mathbb E}[\Phi(\widetilde X_n)]-\widehat{\mathbb E}[\Phi_k(\widetilde X_n)]\right|
    +\left|\widetilde{\mathbb E}[\Phi]-
    \widetilde{\mathbb E}[\Phi_k]
  \right|\\
  \le 
  & \left|
    \widehat{\mathbb E}[\Phi_k(\widetilde X_n)]
    -
    \widetilde{\mathbb E}[\Phi_k]
  \right|+2\|\Phi-\Phi_k\|_K+2 \sup_{n\ge 1}\sup_{P\in\mathcal{P}}P(\widetilde{X}_n\notin K)+2\sup_{\theta\in\mathcal A}P_W(B^\theta\notin K). 
\end{align*}
For any given $\epsilon$, by the tightness, we can choose compact $K$ such that $2\sup_{n\ge1} \sup_{P\in\mathcal{P}}P(\widetilde{X}_n\notin K)+2\sup_{\theta\in\mathcal A}P_W(B^\theta\notin K)<\epsilon/2$. Then we can choose $k$ such that $2\|\Phi-\Phi_k\|_K<\epsilon/2$. Letting $n\to \infty$ yields 
$$ \limsup_{n\to \infty}\left|
    \widehat{\mathbb E}[\Phi(\widetilde X_n)]
    -
    \widetilde{\mathbb E}[\Phi]
  \right|<\epsilon. $$
This proves the theorem.
\end{proof}

\bibliographystyle{plain}
\bibliography{references}

\end{document}